\documentclass[a4,10pt]{article} \usepackage{amssymb}
 \usepackage{amsthm}\usepackage{amsmath}\usepackage{latexsym}\usepackage{srcltx}

\newtheorem{teo}{Theorem}[section]
\newtheorem{oss}[teo]{Remark}
\newtheorem{war}[teo]{Warning}
\newtheorem{Prop}[teo]{Proposition}
\newtheorem{lemma}[teo]{Lemma}

\newtheorem{Defi}[teo]{Definition}

\newtheorem{corollario}[teo]{Corollary}
\newtheorem{no}[teo]{Notation}

\newtheorem{claim}[teo]{Claim}

\newtheorem{claima}[teo]{Problem}

\newcommand{\cc}{_{^{_\HH}}}

\newcommand{\res}{\mathop{\hbox{\vrule height 7pt width .5pt depth 0pt
\vrule height .5pt width 6pt depth 0pt\,}}\nolimits}

\newcommand{\LL}{\mathop{\hbox{\vrule height .5pt width 6pt depth
0pt \vrule height 7pt width .5pt depth 0pt\,}}\nolimits}
\newcommand{\rr}{_{^{_\mathcal{R}}}}

\def \cin{{\mathbf{C}^{\infty}}}

\def\dim {\mathrm{dim}}
\def\dc {{d_{CC}}}

\def\ss{_{^{_{\HS}}}}
\def\eu {_{^{_{\rm Eu}}}}

\def\g{h\cc}

\def\dg{\textit{grad}\cc}
\def\qq{\textit{grad}\ss}

\def \per {\sigma^{2n}\cc}
\def \perh {\sigma^{2n}\cc}
\def \perhu {\sigma^{2}\cc}

\def\Sph{\mathbb{S}_{\mathbb{H}^n}}

\def\UU{\mathcal{U}}

\def \TB{\tau^{_{\TT\!S}}}

\def \nn{\nu\cc}

\def \XX{\mathfrak{X}}

\def \MS{\mathcal{H}\cc}

\def \P{{\mathcal{P}}}

\def \PH{\P\cc}

\def \Om{\Omega}

\def \ee{\mathrm{e}}

\def \R{\mathbb{R}}
\def \Rn{\mathbb{R}^{\DN}}

\def \div{\mathit{div}}

\def \gg{\mathfrak{g}}

\def\divh{\div_{\!^{_\HH}}}

\def\lh{\mathcal{L}\ss}

\def\lg{\mathcal{D}\ss}

\def\tsc{\nabla^{^{_{\TT\!{S}}}}}
\def\gs{\nabla^{_{\HS}}}
\def\gc{\nabla^{_{\HH}}}

\def \Tor{{\textsc{T}}}

\def\UU{\mathcal{U}}

\def\UU{\mathcal{U}}

\def \nn{\nu_{_{\!\HH}}}

\def \ee{\mathrm{e}}

\def \cont{{\mathbf{C}}}
\def \Om{\Omega}
\def \Rn{\mathbb{R}^{n}}
\def \R{\mathbb{R}}

\def \cji {c_{j\,i}(x)}
\def \C { C(x):=[\cji]_{j,i},\,\, {j=1,\ldots,m \,,\, i=1,\ldots,n}}

\def \X {X=(X_{1}, \ldots, X_{m_1})}

\def \X0 {X_{1}(0)\!=\!\partial_{x_{1}}, \ldots, X_{m_1}(0)\!=\!\partial_{x_{m_1}}}

\def \HG {\HH\GG}
\def \HS {\HH\!{S}}

\def \HH {\mathit{H}}

\def \TT {\mathit{T}}
\def \TS {\mathit{T}S}

\def \grad{\textit{grad}}
\def \C0H{\mathbf{C}_{0}^{\infty}(U,\HG)}
\def \C00{\mathbf{C}_{0}^{\infty}(U)}
\def \C01{\mathbf{C}_{0}^{1}(U)}

\def \L1{d\,\mathcal{L}^n}

\def \H1{\mathcal{H}_{{\bf cc}}^{1}}

\def \exp{\textsl{exp\,}}

\def \Om{\Omega}
\def \Rn{\mathbb{R}^{n}}
\def \R{\mathbb{R}}

\def \cji {c_{j\,i}(x)}
\def \C { C(x):=[\cji]_{j,i},\,\, {j=1,\ldots,m \,,\, i=1,\ldots,n}}

\def \gg{\mathfrak{g}}

\def \X {X=(X_{1}, \ldots, X_{m_1})}

\def \X0 {X_{1}(0)\!=\!\partial_{x_{1}}, \ldots, X_{m_1}(0)\!=\!\partial_{x_{m_1}}}

\def \HG {\mathit{H}}

\def \C0H{\mathbf{C}_{0}^{\infty}(\Om,\HG)}
\def \C00{\mathbf{C}_{0}^{\infty}(\Om)}
\def \C01{\mathbf{C}_{0}^{1}(\Om)}

\def \exp{\textsl{exp\,}}

\begin{document}

\title{Hypergeometric solutions of the closed
eigenvalue problem on Heisenberg Isoperimetric Profiles}
\author{Francescopaolo Montefalcone
\thanks{F. M. was partially supported by CIRM, Fondazione Bruno Kessler, Trento, and by the Fondazione CaRiPaRo Project ``Nonlinear Partial Differential Equations: models, analysis, and control-theoretic problems".} }

\maketitle

\date{}

\begin{abstract}
After introducing  the sub-Riemannian geometry of the Heisenberg
group  $\mathbb{H}^n,\,n\geq 1,$ we recall some basics about
hypersurfaces endowed with the $\HH$-perimeter measure $\perh$ and
horizontal Green's formulas. Then, we describe a class of compact
closed hypersurfaces of constant horizontal mean curvature called
\textquotedblleft Isoperimetric Profiles\textquotedblright (they
are not CC-balls!); see Section \ref{Sez3}. Our main purpose is to
study a closed eigenvalue problem on Isoperimetric Profiles, i.e.
$\lh\varphi+\lambda\varphi=0$, where $\lh$ is a 2nd order
horizontal tangential operator analogous to the Laplace-Beltrami
operator; see Section \ref{auts}. This is done starting from the
radial symmetry of Isoperimetric Profiles with respect to a
barycentric axis parallel to the center $T$ of the Lie algebra
$\mathfrak{h}_n$. An interesting feature of radial eigenfunctions
is in that they are hypergeometric functions; see Theorem
\ref{RPF}. Finally, in Section \ref{GENCASF} we shall begin the study of the general case.
\end{abstract}

\tableofcontents

\section{Introduction}\label{basicsHYP}

In the last decades,  sub-Riemannian geometry has become a
subject of great interest due to its
connections with many different areas of Mathematics and Physics,
such as PDE's, Calculus of Variations, Control Theory, Mechanics
and, more recently,  Visual Geometry (see \cite{SCM} and
references therein) and Theoretical Computer Science (see
\cite{Ch3}). For references and comments  we refer the reader to the seminal paper
by Gromov, \cite{Gr1}, Montgomery's book, \cite{Montgomery}, and the
papers by  Pansu, \cite{P4}, Strichartz,
 \cite{Stric}, and Vershik and Gershkovich,
\cite{Ver}.

 In the context of sub-Riemannian geometry, Carnot groups provide a very large class of examples. They are connected, simply connected, nilpotent and stratified\footnote{A vector space $E$  is {\it stratified} if
there exist vector subspaces $E_1,...,E_k$ of $E$ such that
$E=E_1\oplus...\oplus E_k$. The subspace $E_i$ is called the {\it
$i$-th layer} of the stratification while the integer $k$ is the
{\it step} of $E$.} Lie groups which play in sub-Riemannian
geometry a role analogous to that of Euclidean spaces in
Riemannian geometry. A fundamental feature of Carnot groups is
that they are homogeneous groups in the sense of Stein's
definition (see \cite{Stein}), and this means that they are
equipped with an intrinsic family of anisotropic dilations. The
homogeneous dimension of any Carnot group is the weighted sum of
the dimensions of each layer of the stratification and this
integer equals the intrinsic metric dimension of the group,
considered as a metric space with respect to the so-called {\it
Carnot-Carath\'eodory} distance. For an introduction to the
geometry of Lie groups we refer the reader to Helgason,
\cite{Helgason}, and Milnor, \cite{3}, while specifically for
sub-Riemannian geometry to Gromov, \cite{Gr1}, Pansu, \cite{P4},
and Montgomery, \cite{Montgomery}. Carnot groups nowadays
represent an intensive research-field in Analysis and Geometric
Measure Theory; see, for instance, \cite{balogh}, \cite{CDPT},
\cite{vari}, \cite{DanGarN8, gar}, \cite{FSSC3, FSSC5},
\cite{LeoM}, \cite{Monti}, \cite{Mag}, \cite{Monteb},
\cite{Pauls}, \cite{RR}, but these references are far from being
complete.

Among Carnot groups, the most common ones are the so-called {\it
Heisenberg groups} $\mathbb{H}^n$, $n\geq 1$; see
Section \ref{hn}. Roughly speaking, as a manifold, $\mathbb{H}^n$
can be regarded as $\mathbb{C}^n\times \R$  endowed with a suitable polynomial group law $\star$.
Its Lie algebra $\mathfrak{h}_n$ identifies with
$\TT_0\mathbb{H}^n$ (i.e. the tangent space at
$0\in\mathbb{H}^n$). Let us introduce a left invariant frame
$\mathcal{F}=\{X_1, Y_1,...,X_n, Y_n, T\}$ for $\TT\mathbb{H}^n$,
where $X_i=\frac{\partial}{\partial x_i} -
\frac{y_i}{2}\frac{\partial}{\partial t}$,
$Y_i=\frac{\partial}{\partial y_i} +
\frac{x_i}{2}\frac{\partial}{\partial t}$  and
$T=\frac{\partial}{\partial t}$. Denoting by $[\cdot, \cdot]$ the
Lie bracket of vector fields, one has
\[[X_i,Y_i]=T\qquad \mbox{for every}\,\,i=1,...,n,\]and all other commutators
vanish. In other words, $\mathfrak{h}_n$ turns out to be nilpotent
and stratified with center $T$, i.e. $\mathfrak{h}_n=\HH\oplus
\HH_2$, where $\HH={\rm span}_{\R}\{X_1,
Y_1,...,X_i,Y_i,...,X_n,Y_n\}$ is the horizontal bundle and
$\HH_2={\rm span}_{\R}\{T\}$ is the vertical bundle. Later on,
$\mathbb{H}^n$ will be endowed with the  left-invariant Riemannian
metric $h:=\langle\cdot, \cdot\rangle$ making $\mathcal{F}$ an
orthonormal frame. This metric induces a corresponding metric
$h\cc$ on $\HH$ which provides a way to measure  the length of
horizontal curves. In fact, the natural distance in sub-Riemannian
geometry is the {\it Carnot-Carath\'eodory  distance} $\dc$,
defined  by minimizing the (Riemannian) length of all piecewise
smooth horizontal curves joining two different points\footnote{The
last definition is motivated by Chow's Theorem, which implies that
two different points can be joined by (infinitely many) horizontal
curves.}.

The stratification of  $\mathfrak{h}_n$ is connected with the
existence of a 1-parameter group of automorphisms, called {\it
Heisenberg dilations}, defined by $\delta_s (z, t):=(s z, s^2 t)$,
where $(z, t)\in\R^{2n+1}$ denote exponential coordinates of a
point  $p\in\mathbb{H}^n$.  The integer $Q=2n+2$ turns out to be the
(homogeneous) dimension of $\mathbb{H}^n$ as a metric space with
respect to the CC-distance $\dc$.

Let $S\subset\mathbb{H}^n$ be a $\cont^2$-smooth hypersurface. By definition, the
{\it characteristic set} $C_S\varsubsetneq S$ is the set of all points
$x\in S$ such that $\dim \HH_x=\dim (\HH_x \cap \TT_x S)$. In
other words, $x\in S$ is non-characteristic if, and only if, the
horizontal bundle $\HH$ is transversal to $\TS$ at the point $x$.
If $S$ is non-characteristic, the {\it unit $\HH$-normal} along
$S$ is defined by $\nn: =\frac{\PH\nu}{|\PH\nu|}$, where $\nu$ is
the Riemannian unit normal along $S$ and
$\PH:\TT\mathbb{H}^n\longrightarrow\HH$ is the orthogonal
projection onto $\HH$. Moreover, there exists a
homogeneous measure  $\perh$ on $S$, called \it $\HH$-perimeter measure,
\rm which can  be defined in terms of the Riemannian measure
$\sigma\rr^{2n}$ on $S$. More precisely, one has$$\perh \res S =
|\PH \nu |\,\sigma^{2n}\rr\, \res S.$$We have to remark that the measure \it $\perh$ does not see characteristic points.
\rm Furthermore, denoting by $\mathcal{S}_{\bf
cc}^{Q-1}$ the $(Q-1)$-dimensional spherical Hausdorff
measure associated with $\dc$, it turns out that  $\perh(S\cap
B)=k(\nn)\,\mathcal{S}_{\bf cc}^{Q-1}\res
 ({S}\cap B)$ for all  Borel set $B$, where $k(\nn)$ is a density-function called {\it metric factor}; see \cite{Mag}.
On $S\setminus C_S$, there are two important subbundles to be defined: the
{\it horizontal tangent bundle} $\HS\subset \TS$ and the {\it
horizontal normal bundle} $\nn S$. They split the horizontal
bundle $\HH$ into an orthogonal direct sum, i.e. $\HH= \nn S \oplus
\mathit{H}S$. On the other hand, at $C_S\subsetneq S$,  only the subbundle $\HS$ turns out to be defined, and in this case  $\HS=\HH$.

Following the somehow general geometric approach begun in
\cite{Monteb},
 in this paper we address the issue of studying a closed eigenvalue
problem for a
2nd order horizontal tangential operator  $\lh$, which is the
sub-Riemannian analogous to the Laplace-Beltrami operator.
More precisely, let $\qq$ denote the gradient operator on the horizontal tangent bundle $\HS$. A simple way to define
the operator $\lh$ reads as follows:
\[\lh \varphi\,\perh\res S:= d\left( \qq\varphi\LL \perh\right)\res S\]for all $\cont^2$-smooth function $\varphi:S\longrightarrow\R$, where $\qq\varphi=\grad\cc\varphi-\langle\grad\cc\varphi,\nn\rangle\nn$.  In this definition,  $d:\Lambda^{k}(\TT^\ast\mathbb H^n)\longrightarrow\Lambda^{k+1}(\TT^\ast\mathbb H^n)$ is the   exterior derivative and $\LL:\Lambda^{k}(\TT^\ast\mathbb H^n)\longrightarrow\Lambda^{k-1}(\TT^\ast\mathbb H^n)$ is the \textquotedblleft contraction operator\textquotedblright on differential forms; see \cite{Helgason}, \cite{FE}. The importance of the operator $\lh$ comes from some horizontal Green-type formulas discussed in Section \ref{IBPAA}; see also Section \ref{Sez03}.

In a certain sense, what we are trying to do is to prove a sub-Riemannian
generalization of the closed eigenvalue problem for the Laplacian on
Spheres. Actually, we are interested to find explicit solutions to
a \textquotedblleft closed eigenvalue problem\textquotedblright\footnote{This means solving the equation  $\lh\varphi+\lambda\varphi=0$ on a smooth compact (closed) hypersurface $S\subset \mathbb H^n$.} for the operator $\lh$,  in the special case of
{\it Isoperimetric Profiles}. They are compact hypersurfaces
which can be efficiently described in terms of a key-property of
CC-geodesics.

 Let us sketch a picture in the 1st Heisenberg group
$\mathbb{H}^1$. Remind that
 any
CC-geodesic $\gamma$ in $\mathbb{H}^1$ is either a Euclidean
(horizontal) line or a suitable infinite circular helix of
constant slope and whose axis is parallel to the center $T$ of the
Lie algebra $\mathfrak{h}_1$. In the second case, choose a point
${\mathcal S} \in \gamma$ and take the (positively oriented)
$T$-line $\mathcal V$ over this point. There exists a first
consecutive point to $\mathcal S$, denoted by $\mathcal N$,
which\footnote{This point can be interpreted as the cut point of
$\mathcal S$ along $\gamma$. In fact this is the end-point of all
CC-geodesics starting from $\mathcal S$ with same slope. Note
however that, strictly speaking, the cut locus of any point in
$\mathbb{H}^n\, (n\geq 1)$ coincides with the vertical $T$-line
over the point.} belongs to $\gamma \cap \mathcal V$. These
points, henceforth called South and North poles, determine a
minimizing connected subset of $\gamma$.  The slope of this helix
is uniquely determined by the CC-distance of the  poles. Now
rotating (around the vertical axis from $\mathcal{N}$ to
$\mathcal{S}$) the curve joining the poles yields a closed convex
surface very similar to an ellipsoid, henceforth called {\it
Isoperimetric Profile}. A similar description can be done  in the
general case. Hereafter, we shall denote by $\Sph$ the unit
Isoperimetric Profile of $\mathbb H^n$; see Section \ref{Sez3}.

 We also  stress that Isoperimetric Profiles  are
constant horizontal mean curvature playing in Heisenberg groups an equivalent role of spheres in Euclidean spaces. In this regard,   there is a longstanding
conjecture claiming that Isoperimetric Profiles
minimize the $\HH$-perimeter measure among \textquotedblleft finite
$\HH$-perimeter sets\textquotedblright \footnote{In the variational sense; for a definition, see \cite{FSSC3}.} having fixed volume; see \cite{CDPT},
\cite{DanGarN8}, \cite{gar}, \cite{gar2}, \cite{LeoM},
\cite{Monti, Monti2}, \cite{Monti3}, \cite{Monti4}, \cite{Ni},
\cite{P1, Pansu2}, \cite{RR}. For this reason, it seems interesting to study some features of these sets from a sub-Riemannian point of view.

The plan of the paper is the following.
 After introducing in Section \ref{hn}
 the sub-Riemannian geometry of
Heisenberg groups $\mathbb{H}^n,\,n\geq 1,$ we overview some
basic facts concerning the theory of hypersurfaces endowed with the $\HH$-perimeter measure $\perh$.
We also prove some useful geometric identities for the sequel; see
Section \ref{sez2}. In Section \ref{IBPAA} and Section \ref{Sez03} we discuss some
horizontal divergence formulas.  In Section \ref{Sez3} we briefly describe Heisenberg
Isoperimetric Profiles. In Section
\ref{auts} we examine some general features of the closed eigenvalue problem  for the operator $\lh$ on
compact closed hypersurfaces, i.e.
\begin{displaymath}
\left\{%
\begin{array}{ll}
    \quad\,\lh\varphi=-\lambda\,\varphi \quad\mbox{on}\,\,S\setminus C_S\\
    \int_S\varphi\,\perh=0;\\
\end{array}%
\right..\end{displaymath}see Problem \ref{PEH}. As already said, our main interest concerns the case of Isoperimetric Profiles and so we shall choose $S=\Sph$.

\begin{oss}Let  $C_P$  denote the best constant in the following Poincar\'e-type inequality
$$\int_{\Sph}\varphi^2 \perh\leq C_P\int_{\Sph} |\qq \varphi|^2 \perh$$ for all $\varphi$ belonging to the horizontal Sobolev space $\mathcal{H}(\Sph, \perh)$ such that $\int_{\Sph}\varphi\,\perh=0$. Then, it turns out that $$C_P=\frac{1}{\mu},$$where $\mu$ is the 1st eigenvalue  of Problem \ref{PEH}; see Section \ref{auts}.
\end{oss}

Our starting point is the radial symmetry of Isoperimetric
Profiles with respect to a barycentric axis parallel to the center
$T$ of $\mathfrak{h}_n$; see Section
\ref{SEZPr}. In this way, we reduce ourselves to a \it radial
 eigenvalue problem; \rm see Problem \ref{prad}.
Therefore, we have to study an O.D.E. of hypergeometric
type, which can be explicitly integrated by classical methods
for hypergeometric equations. The spectrum
of $\lh$ on radial functions is found by using some natural  mixed conditions for this problem. Eigenfunctions turn out to be
hypergeometric functions; see Theorem \ref{RPF}. Then, after studying the radial case, in
Section \ref{GENCASF} we overview some features of the general
case. Among other things, we prove
that the spectrum of $\lh$  contains
the radial spectrum; see Proposition \ref{teoauto}. This is done
by showing that the spherical mean of each eigenfunction of $\lh$
 associated with an eigenvalue $\lambda$ (whenever different from $0$) must be an eigenfunction of
the radial eigenvalue problem associated with $\lambda$. This fact raises the following question (see e.g. Remark \ref{zxc3}):
\begin{center}\it Is the 1st eigenvalue of Problem \ref{PEH}  equal to the first eigenvalue of  Problem \ref{prad}? \end{center}
Finally, by applying a standard trick (see \cite{FCS}) we will prove some related  inequalities; see Corollary \ref{zxc2} and Corollary \ref{zxc3}. These inequalities are stronger than what one might expect and this fact seems to suggest a positive answer to the previous question.\\\\

\begin{center}\sc Acknowledgements \end{center}
{\small This paper is the first  of a series devoted to study, in the
context of Heisenberg groups, constant  horizontal mean
curvature hypersurfaces and among them, Isoperimetric Profiles.
In many respects, this project begun during my visit at Purdue
University in November 2006.
So I wish to express my gratitude to Prof. N. Garofalo for his kind invitation and for many interesting conversations during the last years. Furthermore,
I would like  to thank Prof. A. Parmeggiani for his sincere
interest in this work and C. Senni for his help with \sc Mathematica.}

\subsection{Heisenberg group $\mathbb{H}^n$}\label{hn}

The {\it $n$-th Heisenberg group} $(\mathbb{H}^n,\star)$, $n\geq
1$, is a connected, simply connected, nilpotent and stratified Lie
group of step 2 on $\R^{2n+1}$, with respect to a polynomial
group law $\star$; see below. Its {\it Lie algebra}
$\mathfrak{h}_n$  is isomorphic to the $(2n+1)$-dimensional real vector space $\R^{2n+1}$ and it will be identified with the tangent space at
the identity $0$  of $\mathbb{H}^n$, i.e.
$\mathfrak{h}_n\cong\TT_0\mathbb{H}^n$. We shall
adopt {\it exponential coordinates of the 1st kind}. Hence every
point $p\in\mathbb{H}^n$ identifies with an ordered
$(2n+1)$-tuple of the Lie algebra $\mathfrak{h}_n$, i.e.
$p=\exp(x_1,y_1,...,x_i,y_i,...,x_n,y_n, t)$. For simplicity, $(z, t)\in\R^{2n+1}$ will denote exponential coordinates of a
generic point $p\in\mathbb H^n$.  In order to
describe the algebraic structure of $\mathfrak{h}_n$, let us fix
a global left-invariant frame $\mathcal{F}:=\{X_1,Y_1,...,X_i,Y_i,...,X_n,Y_n,T\}$ for
$\TT\mathbb{H}^n$, where we have set
\begin{eqnarray*}X_i(p):=\frac{\partial}{\partial x_i} -
\frac{y_i}{2}\frac{\partial}{\partial t},\qquad
 Y_i(p):=\frac{\partial}{\partial y_i} +
\frac{x_i}{2}\frac{\partial}{\partial t}\quad
(i=1,...,n),\quad T(p) := \frac{\partial}{\partial t},\end{eqnarray*}for
every $p\in\mathbb{H}^n$. Denoting by $[\cdot, \cdot]$ the Lie
bracket of vector fields, one has
$[X_i,Y_i]=T$ for every $i=1,...,n$, and all other commutators
vanish. Hence,
 $\mathfrak{h}_n$ turns out to be a nilpotent and stratified Lie algebra of step 2 with center $T$, i.e.
$\mathfrak{h}_n=\HH\oplus \HH_2$, where
$\HH:={\rm span}_{\R}\{X_1, Y_1,...,X_i,Y_i,...,X_n,Y_n\}
$ and $\HH_2={\rm span}_{\R}\{T\}$. The first layer $\HH$ is called
{\it horizontal}, whereas the second layer, spanned by
$T$, is called {\it vertical}; they are both smooth
subbundles of $\TT\mathbb{H}^n$.

The {\it Baker-Campbell-Hausdorff formula} uniquely determines the
group law $\star$ of $\mathbb{H}^n$, starting from the
``structure'' of its Lie algebra $\mathfrak{h}_n$; see
\cite{Corvin}. More precisely, we have
$$\exp X\star\exp Y=\exp(X\diamond Y)\quad\mbox{for every}\,\,X,\,Y \in \mathfrak{h}_n,$$ where
$\diamond:\mathfrak{h}_n\times \mathfrak{h}_n\longrightarrow
\mathfrak{h}_n$ is the operation defined by $X\diamond Y:= X + Y+
\frac{1}{2}[X,Y]$. Therefore, for every
 $p=\exp(x_1,y_1,...,x_n,y_n,
t),\,\,p'=\exp(x'_1,y'_1,...,x'_n,y'_n, t')\in \mathbb{H}^n$, it follows that
\begin{equation}\label{1}
 p\star p':=\exp\left(x_1+x_1', y_1+y_1',...,x_n+x_n', y_n+y_n', t+t'+ \frac{1}{2}\sum_{i=1}^n
\left(x_i y'_{i}- x'_{i} y_i\right)\right).
\end{equation}Hence, the inverse of any
${p}\in\mathbb{H}^n$  is ${p}^{-1}:=\exp(-{x}_1,-y_1...,-{x}_{n},
-y_n, -t)$ and $0=\exp(0_{\R^{2n+1}})$.

\begin{Defi}\label{dccar} A {\rm sub-Riemannian metric} $\g$  is
a symmetric positive bilinear form on the horizontal bundle
$\HH\subset\TT\mathbb{H}^n$. The {\rm {CC}-distance} $\dc(p,p')$
between $p, p'\in \mathbb{H}^n$ is defined by
$$\dc(p, p'):=\inf \int\sqrt{\g(\dot{\gamma},\dot{\gamma})} dt,$$
where the infimum is taken over all piecewise-smooth horizontal curves
$\gamma$ joining $p$ to $p'$.
\end{Defi}

By Chow's Theorem, every couple of points can be connected by a
horizontal curve (not unique, in general) and this implies
that $\dc$ is a metric on $\mathbb{H}^n$. The  topology generated by the CC-metric turns out to be
equivalent to the standard  topology on $\R^{2n+1}$; see
\cite{Gr1}, \cite{Montgomery}.

From now on, we will equip $\TT\mathbb{H}^n$ with the left-invariant
Riemannian metric
  $h:=\langle\cdot,\cdot\rangle$ making $\mathcal{F}$ an
orthonormal frame and  assume $\g:=h|_\HH.$

The {\it structural constants}\footnote{Let $\gg\cong\Rn$ be a Lie
algebra  with respect to $[\cdot, \cdot]$, let $\langle\cdot,
\cdot\rangle$ be a Euclidean metric on $\gg$ and let
$\mathcal{F}=\{X_1,...,X_n\}$ be an o.n.  basis of $\gg$. The
{\it structural constants} of $\gg$ associated with $\mathcal{F}$
are $n^3$ smooth functions defined by
$$C^r_{ij}:=\langle [X_i,X_j],
 X_r\rangle\qquad i, j, r=1,...,n.$$
The structural constants embody all algebraic features of $\gg$;
see \cite{Helgason}. They satisfy the following properties:
\begin{itemize}\item[{\rm (i)}] $C^r_{ij} +C^r_{ji}=0$\quad (Skew-symmetry) \item[{\rm (ii)}]
 $\sum_{j=1}^{n} C^i_{jl}C^{j}_{rm} +
C^i_{jm}C^{j}_{lr} + C^i_{jr}C^{j}_{ml}=0$\quad (Jacobi
Identity).\end{itemize}} of the Heisenberg Lie algebra
$\mathfrak{h}_n$ are completely described by  the
skew-symmetric $(2n\times 2n)$-matrix
$$C\cc^{2n+1}:=\left|
\begin{array}{cccccccc}
  0 & 1 & 0 & 0 & \cdot &  \cdot &  0 &  0 \\
  -1 & 0 & 0 & 0 &  \cdot &  \cdot &  0 &  0 \\
  0 & 0 & 0 & 1 &  \cdot & \cdot & 0 & 0 \\
  0 & 0 & -1 & 0 & \cdot & \cdot & 0 & 0 \\
  \cdot & \cdot & \cdot & \cdot & \cdot & \cdot & \cdot & \cdot \\
  \cdot & \cdot & \cdot & \cdot & \cdot & \cdot & \cdot & \cdot \\
  0 & 0 & 0 & 0 & \cdot & \cdot & 0 & 1 \\
  0 & 0 & 0 & 0 & \cdot & \cdot & -1 & 0
\end{array}%
\right|.$$We remark that $C\cc^{2n+1}$ is the matrix associated
with the skew-symmetric bilinear map
$\Gamma\cc:\HH\times\HH\longrightarrow \R$ defined by
$\Gamma\cc(X, Y)=\langle[X, Y], T\rangle$.

 For
every $p\in\mathbb{H}^n$, we denote by
$L_p:\mathbb{H}^n\longrightarrow\mathbb{H}^n$ the {\it
left translation by $p$}, i.e. $L_pp':=p\star p'$. The map $L_p$
is a group homomorphism and its differential
${L_p}_\ast:\TT_0\mathbb{H}^n\longrightarrow\TT_p\mathbb{H}^n$  is
given by

\begin{displaymath}
{L_p}_\ast=\frac{\partial (p\star p')}{\partial p'}\bigg|_{p'=0}=
\left[%
\begin{array}{ccccccccc}
  1 & 0 & \ldots &0&0& \ldots& 0 & 0 & 0 \\
  0 & 1 & \ldots &0&0& \ldots & 0 &0 & 0 \\
  \vdots & \vdots & \ldots & \vdots & \vdots & \ldots & \vdots & \vdots&\vdots \\
  0 & 0 & \ldots & 0& 0 & \ldots & 0 & 1 & 0 \\
 - \frac{y_1}{2} & +\frac{x_1}{2} & \ldots&-\frac{y_i}{2} & +\frac{x_i}{2} & \ldots & -\frac{y_n}{2} & +\frac{x_n}{2} & 1\\
\end{array}%
\right].
\end{displaymath}Equivalently, one has ${L_p}_\ast={\rm{col}}[X_1(p),Y_1(p),...X_n(p),Y_n(p),T(p)]$.
A key feature of $\mathbb{H}^n$ is
that there exists a 1-parameter group of automorphisms
$\delta_s:\mathbb{H}^n \longrightarrow\mathbb{H}^n\,(s\geq 0)$,
hereafter called {\it Heisenberg dilations}, defined by
$\delta_s p :=\exp\left(s z, s^2 t\right)$ for every $s\geq
0$, where $p=\exp(z, t)$.  As already said, the {\it
homogeneous dimension} of $\mathbb{H}^n$ with respect to the intrinsic dilations is the integer given
by $Q:=2n+2$ which
equals the {\it Hausdorff dimension} of $\mathbb{H}^n$ as a metric
space with respect to the CC-distance $\dc$; see also
\cite{Montgomery}.

We shall denote by $\nabla$ the unique {\it left-invariant
Levi-Civita connection} on $\mathbb{H}^n$ associated with the
metric $h=\langle\cdot,\cdot\rangle$. We recall that, for
every $X, Y, Z\in \XX:=\cin(\mathbb{H}^n, \TT\mathbb{H}^n)$ one
has
\[\left\langle\nabla_XY,Z\right\rangle=\frac{1}{2}
\left(\langle[X, Y], Z\rangle-\langle[Y, Z], X\rangle + \langle[Z,
X], Y\rangle\right);\]for more details, see \cite{3}. For every $X, Y\in\XX\cc:=\cin(\mathbb{H}^n,\HH)$, we shall set
$\gc_X Y:=\PH\left(\nabla_X Y\right)$, where $\PH$ denotes the orthogonal
projection onto $\HH$.  The operation $\gc$ is a
vector-bundle connection, called {\it $\HH$-connection}; see
\cite{Monteb} and references therein. It is not difficult to show
that $\gc$ is {\it flat} and {\it compatible with the
sub-Riemannian metric} $\g$, i.e.
$X\langle Y, Z \rangle=\left\langle \gc_X Y, Z \right\rangle +
\left\langle Y, \gc_X Z \right\rangle$ for all $X,
Y\in \XX\cc$. Furthermore, $\gc$ is {\it torsion-free}, i.e.
 $\gc_X Y - \gc_Y X-\PH[X,Y]=0$ for all $X, Y\in \XX\cc$. All these properties  follow from the very definition of
$\gc$ by using standard features of the Levi-Civita connection
$\nabla$ on $\mathbb{H}^n$.

\begin{Defi}
For every $\psi\in\cin(\mathbb{H}^n)$  the {\rm $\HH$-gradient} of
$\psi$ is the unique horizontal vector field $\dg \psi\in\XX\cc$
such that
$\langle\dg \psi,X \rangle= d \psi (X) = X
\psi$ {for all} $X\in \HH.$ The {\rm $\HH$-divergence}
$\divh X$ of any $X\in\XX\cc$ is given, at each point $p\in
\mathbb{H}^n$, by
$\divh X(p):= \mathrm{Trace}\left(Y\longrightarrow \gc_{Y} X
\right)(p)$ $(Y\in \HH_p)$. The {\rm $\HH$-Laplacian}
$\Delta\cc$  is the 2nd order differential operator defined by
\[
\Delta\cc\psi := \div\cc(\dg\psi)\quad\mbox{for every}\,\,\psi\in
\cin(\mathbb{H}^n).\]
\end{Defi}

We now recall some basic notions about left invariant forms on
$\mathbb{H}^n$.  Starting from the left invariant frame
$\mathcal{F}$ and by duality\footnote{The duality is understood
with respect to the fixed left-invariant metric $h$ on
$\mathbb{H}^n$. More generally, if $(M, h)$ is a Riemannian
manifold and $X\in\TT M$, then $X^\ast(Y):=h(X, Y)$ for every
$Y\in\TT M$.}, we define  a global coframe
$\mathcal{F}^\ast:=\{X^\ast_1,Y^\ast_1,...,X^\ast_i,Y^\ast_i,...,X^\ast_n,Y^\ast_n,T^\ast\}$
of {\it left invariant $1$-forms} for the cotangent bundle
$\TT^\ast\mathbb{H}^n$, where
\begin{eqnarray*}X_i^\ast= dx_i,\qquad  Y_i^\ast=dy_i\quad (i=1,...,n),\qquad T^\ast=dt
+ \frac{1}{2}\sum_{i=1}^n\left(y_i d x_i  -  x_i d
y_i\right).\end{eqnarray*}We shall set
$\theta:=T^\ast$ to denote the so-called {\it contact $1$-form} of
$\mathbb{H}^n$.

Finally the top-dimensional {\it left-invariant volume form}
$\sigma^{2n+1}\rr$ of $\mathbb{H}^{n}$ is defined by
$$\sigma^{2n+1}\rr:=\left(\bigwedge_{i=1}^n dx_i \wedge d y_i\right) \wedge
\theta.$$The measure obtained by integration of $\sigma^{2n+1}\rr$
equals the {\it Haar measure} of $\mathbb{H}^{n}$.

\subsection{Hypersurfaces, measures and some useful formulas}\label{sez2}

Let $S\subset\mathbb{H}^n$ be a $\cont^1$-smooth hypersurface and
denote by $\nu$ its (Riemannian) unit normal. We have
$\nu=\sum_{i=1}^n \big(\langle\nu, X_i\rangle X_i + \langle\nu,
Y_i\rangle Y_i \big) + \nu_T\, T$, where
\begin{eqnarray*}\langle\nu, X_i\rangle&=&\frac{\langle X_i, {\rm n}\eu\rangle_{\R^{2n+1}}}
{\sqrt{\sum_{i=1}^n\left(\langle X_i, {\rm
n}\eu\rangle_{\R^{2n+1}}^2 + \langle Y_i, {\rm
n}\eu\rangle_{\R^{2n+1}}^2\right)+ \langle T, {\rm
n}\eu\rangle_{\R^{2n+1}}^2}}\quad i=1,...,n,\\
\langle\nu, Y_i\rangle&=&\frac{\langle Y_i, {\rm
n}\eu\rangle_{\R^{2n+1}}}{\sqrt{\sum_{i=1}^n\left(\langle X_i,
{\rm n}\eu\rangle_{\R^{2n+1}}^2 + \langle Y_i, {\rm
n}\eu\rangle_{\R^{2n+1}}^2\right)+ \langle T, {\rm
n}\eu\rangle_{\R^{2n+1}}^2}}\quad i=1,...,n,\\
 \langle\nu, T\rangle&=&\frac{\langle T, {\rm
n}\eu\rangle_{\R^{2n+1}}}{\sqrt{\sum_{i=1}^n\left(\langle X_i,
{\rm n}\eu\rangle_{\R^{2n+1}}^2 + \langle Y_i, {\rm
n}\eu\rangle_{\R^{2n+1}}^2\right)+ \langle T, {\rm
n}\eu\rangle_{\R^{2n+1}}^2}}=:\nu_{T}.\end{eqnarray*}In the above formulas,
$\rm{n}\eu$ denotes the Euclidean unit normal along $S$
and $\langle \cdot, \cdot\rangle_{\R^{2n+1}}$ is the
standard metric on $\R^{2n+1}$.

The Riemannian measure $\sigma^{2n}\rr$ on a smooth
hypersurface $S$ can be defined as {\it contraction}\footnote{Let $M$
be a Riemannian manifold. The linear map $\LL: \Lambda^r(\TT^\ast
M)\rightarrow\Lambda^{r-1}(\TT^\ast M)$ is defined, for $X\in \TT
M$ and $\omega^r\in\Lambda^r(\TT^\ast M)$, by $(X \LL \omega^r)
(Y_1,...,Y_{r-1}):=\omega^r (X,Y_1,...,Y_{r-1})$; see, for
instance, \cite{FE}. This operation is called {\it contraction} or
{\it interior product}.} of the top-dimensional volume form
$\sigma^{2n+1}\rr$ by the unit normal $\nu$ along $S$, i.e.
$\sigma^{2n}\rr\res S := (\nu\LL \sigma^{2n+1}\rr)|_S$.

We say that $p\in S$ is a {\it characteristic point}  if
$\dim\,\HH_p = \dim (\HH_p \cap \TT_p S)$. By definition, the {\it
characteristic set} of $S$ is the set of all characteristic
points, i.e.
$C_S:=\{x\in S : \dim\,\HH_p = \dim (\HH_p \cap \TT_p S)\}.$ It is worth noting that $p\in C_S$ if,
 and only if, $|\PH\nu(p)|=0$. Since $|\PH\nu(p)|$ is continuous
 along $S$, it follows that
$C_S$ is a closed subset of $S$, in the relative topology. We  also remark that characteristic points are few, since under the
preceding assumptions, the $(Q-1)$-dimensional Hausdorff measure
of $C_S$ vanishes, i.e. $\mathcal{H}_{CC}^{Q-1}(C_S)=0$; see, for
instance, \cite{balogh}, \cite{Mag}.
\begin{oss}\label{CSET}Let $S\subset \mathbb{H}^n$ be a $\cont^2$-smooth hypersurface.
A straightforward application of {\rm Frobenius Theorem}
 about integrable distributions shows  that the topological dimension
 of $C_S$ is strictly less than $(n+1)$; see also \cite{Gr1}. For other results  about the size of $C_S$
  in $\mathbb{H}^n$, see \cite{balogh}. \end{oss}

In the sequel we will need another measure on hypersurfaces, the
so-called {\it $\HH$-perimeter measure}; see \cite{FSSC3},
\cite{DanGarN8, gar}, \cite{Mag}, \cite{Monte, Monteb},
\cite{P1},  \cite{Pauls}, \cite{RR} and references therein.

\begin{Defi}[$\perh$-measure]\label{sh}
Let $S\subset\mathbb{H}^n$ be a $\mathbf{C}^1$-smooth
non-characteristic
 hypersurface and let $\nu$ be the unit
normal vector along $S$. The normalized projection of $\nu$ onto
$\HH$ is called {\rm unit $\HH$-normal} along $S$, i.e. $\nn:
=\frac{\PH\nu}{|\PH\nu|}$.  The {\rm $\HH$-perimeter form}
$\perh\in\bigwedge^{2n}(\TT^\ast S)$ is the contraction of the
volume form $\sigma^{2n+1}\rr$ of $\mathbb{H}^n$ by the
horizontal unit normal $\nn$, i.e. $\perh \res S:=\left(\nn \LL
\sigma^{2n+1}\rr\right)\big|_S.$\end{Defi} If
$C_S\neq\emptyset$ we  extend $\perh$ up to $C_S$  by
setting $\perh\res C_{S}= 0$. Moreover, it turns out that$$\perh \res S = |\PH \nu |\,\sigma^{2n}\rr\, \res S.$$
At each point $p\in S\setminus {C}_S$ one has $\HH_p= {\rm
span}_\R\{\nn(p)\} \oplus \mathit{H}_p S$, where  $\mathit{H}_p
S:=\HH_p\cap\TT_p S$. So we can define, in the obvious way, the
associated subbundles $\HS \subset \TS$ and $\nn S$  called,
respectively, {\it horizontal tangent bundle} and {\it horizontal
normal bundle} along $S$.

\begin{oss}\label{SERES} We have $\dim \HH_p S=\dim \HH-1=2n-1$ at each point $p\in S\setminus C_S$. Note that the  definition of $\HS$ makes sense even if $p\in C_S$, but in such a case $\dim \HH_p S=\dim \HH_p =2n$.
\end{oss}

With respect to the horizontal o.n. frame $\mathcal{F}\cc:=\{X_1,
Y_1,..., X_n, Y_n\}$, the unit $\HH$-normal $\nn$ can
be written out as $\nn=\sum_{i=1}^n\big(\langle\nn, X_i\rangle X_i
+ \langle\nn, X_i\rangle Y_i\big)$, where

\begin{eqnarray*}\langle\nn, X_i\rangle&:=&\frac{\langle X_i, {\rm
n}\eu\rangle_{\R^{2n+1}}}{\sqrt{\sum_{i=1}^n\left({\langle X_i,
{\rm n}\eu\rangle^2_{\R^{2n+1}}} + {\langle Y, {\rm
n}\eu\rangle^2_{\R^{2n+1}}}\right)}}\quad i=1,...,n,\\\langle\nn,
Y_i\rangle&:=&\frac{\langle Y_i, {\rm
n}\eu\rangle_{\R^{2n+1}}}{\sqrt{\sum_{i=1}^n\left({\langle X_i,
{\rm n}\eu\rangle^2_{\R^{2n+1}}} + {\langle Y, {\rm
n}\eu\rangle^2_{\R^{2n+1}}}\right)}}\quad
i=1,...,n.\end{eqnarray*}Another important geometric object (see \cite{Monte, Monteb}, \cite{gar}) is given by
$\varpi:=\frac{\nu_{T}}{|\PH\nu|}$. Clearly,  $\varpi$ is not defined at $C_S$, but one has $\varpi\in L^1_{loc}(S,
\perh)$.
\begin{no} \label{notgepperp}We set
$z^\perp:=-C^{2n+1}\cc z=(-y_1, x_1, ...,-y_n,
x_n)\in\R^{2n}$  and $X^\perp:=-C^{2n+1}\cc X$, for every $X\in\HH$. In particular, the horizontal tangent vector field $\nn^\perp$ is called \rm characteristic direction \it  along $S$. We also set
$C\cc(\varpi):=\varpi\,C\cc^{2n+1}$. Note that
$C\cc(\varpi)\nn=-\varpi\nn^\perp=-\varpi \sum_{i=1}^n\left(-\langle\nn, Y_i\rangle X_i +
\langle\nn, X_i\rangle Y_i\right)$.
\end{no}
\begin{Defi}\label{movadafr}Let $S\subset\mathbb{H}^n$ be a $\cont^2$-smooth non-characteristic hypersurface.
We call {\rm adapted  frame} along $S$ any o.n. frame
${\mathcal{F}}:=\{\tau_1,...,\tau_{2n+1}\}$ for
$\TT\mathbb{H}^n$
 such that:\begin{center}{\rm (i)} $\tau_1|_S=\nn$,\quad {\rm (ii)} $\HH_pS=\mathrm{span}_\R\{\tau_2(p),...,\tau_{2n}(p)\}$ for every $p\in S$,\quad {\rm (iii)} $\tau_{2n+1}:= T$.
\end{center}We
also set $I\cc:=\{1,2,3,...,2n\}$ and
 $I\ss:=\{2,3,...,2n\}$.\end{Defi}

\begin{lemma}\label{Sple} Let $S\subset\mathbb{H}^n$  be a $\cont^2$-smooth
non-characteristic hypersurface and $p\in S$. Then, we can always
choose an adapted o.n. frame
${\mathcal{F}}=\{\tau_1,...,\tau_{2n+1}\}$ along $S$ such that
$\langle\nabla_{X}{\tau_i}, \tau_j\rangle=0$ at $p$ for every $i,
j\in I\ss$ and  every $X\in\HH_{p}S$.
\end{lemma}
For a proof, see  Lemma 3.8 in \cite{Monteb}.
The next definitions can  be found in \cite{Monteb}, for
$k$-step Carnot groups. {\it Later on, unless otherwise specified,
we shall assume that $S\subset\mathbb{H}^n$ is a $\cont^2$-smooth
non-characteristic hypersurface}.

\begin{Defi}Let $i=1, 2$. We shall denote by
$\cont^i\ss(S)$ the space of functions whose $i$-th
$\HS$-derivatives are continuous on $S$.  Analogously, for any open subset $\UU\subseteq S$, we set $\cont^i\ss(\UU)$, to denote the space of functions whose $i$-th
$\HS$-derivatives are continuous in $\UU$. In particular, since $C_S$ is closed, we may take $\UU=S\setminus C_S$.
\end{Defi}

\begin{war}The last definition extends to the case $C_S\neq \emptyset$ by requiring that all  $i$-th
$\HS$-derivatives be continuous at each characteristic point $p\in C_S$; see also Remark \ref{SERES}.
 \end{war}

 Let $\tsc$ denote the connection
on $S$ induced from the Levi-Civita connection $\nabla$ on
$\mathbb{H}^n$. We define a partial connection $\gs$
associated with the subbundle $\HS\subset\TT S$ by
setting\footnote{$\P\ss:\TT{S}\longrightarrow\HS$ denotes the
orthogonal projection operator of $\TT{S}$ onto $\HS$.}
$$\gs_XY:=\P\ss\left(\tsc_XY\right)\quad\mbox{for every}\,\,X,Y\in\XX^1\ss:=\cont^1(S, \HS)$$
 Starting from the orthogonal decomposition
 $\HH=\nn S\oplus\HS$,
it can be shown that$$\gs_XY=\gc_X Y-\left\langle\gc_X
Y,\nn\right\rangle \nn\quad\mbox{for every}\,\,X,
Y\in\XX^1\ss.$$\begin{Defi}\label{H2ffHmc}Given $\psi\in
\cont^1\ss(S)$, we define the $\HS$-{\rm gradient of $\psi$}   to be the unique
horizontal tangent vector field $\qq\psi\in\XX^0\ss:=\cont(S, \HS)$ such that
$\langle\qq\psi,X \rangle= d \psi (X) = X
\psi$ {for every} $X\in \HS$. The {\rm
$\HS$-divergence} $\div\ss X$ of any  $X\in\XX^1\ss$ is defined, at each
point $p\in S$, by
$$\div\ss X (p) := \mathrm{Trace}\left(Y\longrightarrow
\gs_Y X \right)(p)\quad\,(Y\in \HH_pS)$$and one has $\div\ss X\in\cont(S)$. The {\rm $\HS$-Laplacian}
$\Delta_{_{\HS}}$ is the 2nd order differential operator given by
\[\Delta\ss\psi := \div\ss(\qq\psi)\quad\mbox{for every}\,\,\psi\in
\cont^2\ss(S).\]
\end{Defi}

\begin{Defi}\label{curvmed}
The {\rm horizontal 2nd fundamental form} of $S$, ${B\cc}:
\XX^1\ss\times\XX^1\ss\longrightarrow\cont(S)$,  is the bilinear
map given by
$B\cc(X,Y):=\left\langle\gc_X
Y,\nn\right\rangle$ {for every} $X,\,Y\in\XX^1\ss$. The {\rm horizontal mean
curvature} is the trace of $B\cc$, i.e. $\MS:=\mathrm{Tr}{B\cc}$.
The {\rm torsion} $\Tor\ss$ of  $\gs$ is defined by
$\Tor\ss(X,Y):=\gs_XY-\gs_YX-\PH[X,Y]$ {for every} $X,\,Y\in\XX^1\ss$.\end{Defi}

If $n=1$ the torsion is zero. Nevertheless, if $n>1$ the torsion
does not vanish, in general, because  $B\cc$ is {\it not
symmetric}; see \cite{Monteb}.
\begin{no}
Let $X\in\XX\cc$  and let $\phi\in\cont^2(\mathbb H^n)$.  We shall denote by $X(\phi):=\langle \grad\cc \phi, X\rangle$ and
 $X^{(2)}(\phi):=X(X(\phi))$, the 1st and 2nd derivatives of $\phi$ along $X$, respectively. Furthermore, functions defined on $S$ will be thought of as restrictions $\phi|_S$ to $S$ of functions defined on $\mathbb H^n$.
\end{no}

Now we prove some identities useful for the sequel.
\begin{lemma}\label{ljjjkl} Let $S\subset\mathbb{H}^n$ be a $\cont^2$-smooth
non-characteristic hypersurface and let
 $\phi\in\cin(\mathbb H^n)$. Then
\begin{eqnarray}\label{ljjjkltay}\Delta\ss\phi=\Delta\cc\phi +
\MS\frac{\partial\phi}{\partial\nn}-\left\langle {\rm
Hess}\cc \phi\, \nn, \nn\right\rangle.\end{eqnarray}
\end{lemma}
\begin{proof}Using an adapted frame ${\mathcal{F}}$, we compute
 \begin{eqnarray*}\Delta\cc\phi&=&\sum_{i\in I\cc}\left(\tau^{(2)}_i-\gc_{\tau_i}\tau_i\right)(\phi)
 \\&=&\tau^{(2)}_1
 (\phi)-\left(\gc_{\tau_1}\tau_1\right)(\phi)+\sum_{i\in I\ss}\left(\left(\tau^{(2)}_i-
 \gs_{\tau_i}\tau_i\right)(\phi)-\left\langle\gc_{\tau_i}\tau_i,
 \tau_1\right\rangle\tau_1(\phi)\right)\\&=&\tau^{(2)}_1
 (\phi)-\left(\gc_{\tau_1}\tau_1\right)(\phi)+ \Delta\ss\phi
 - \MS\tau_1(\phi).
 \end{eqnarray*}Note that the first and the last identities come from the usual invariant definition of the Laplace
 operator on Riemannian manifolds (or vector bundles); see \cite{Hicks}.
 Now we claim that \[\tau^{(2)}_1
 (\phi)-\left(\gc_{\tau_1}\tau_1\right)(\phi)=\left\langle{\rm
Hess}\cc \phi\,\tau_1, \tau_1\right\rangle.\]Assuming $\tau_1=\sum_{i\in
I\cc}A^1_iX_i$  yields
\begin{eqnarray*}\tau^{(2)}_1(\phi)=\sum_{i\in
I\cc}\tau_1(A^1_iX_i(\phi))=\sum_{i,j\in
I\cc}\left(\tau_1(A^1_i)X_i(\phi) +
A^1_iA^1_jX_j(X_i(\phi))\right).\end{eqnarray*}Since
\[\gc_{\tau_1}\tau_1=\sum_{i, j\in I\cc}\left(\tau_1(A^1_i)X_i+
A^1_iA^1_j\underbrace{\gc_{X_i}X_j}_{=0}\right),\]the claim
follows because
\begin{eqnarray}\label{idder}\tau^{(2)}_1
 (\phi)-\left(\gc_{\tau_1}\tau_1\right)(\phi)=\sum_{i,
j\in I\cc} A^1_iA^1_jX_j(X_i(\phi))=\langle {\rm Hess}\cc
\phi\,\tau_1,\tau_1\rangle.\end{eqnarray}
 \end{proof}

\begin{lemma}\label{lppl}Let $S\subset\mathbb{H}^n$ be a $\cont^2$-smooth
non-characteristic hypersurface. Then
$$\gc_{\nn}\nn=-\frac{\qq\varpi}{\varpi}+\varpi\nn^\perp.$$
\end{lemma}
\begin{proof}We make use of an adapted frame ${\mathcal{F}}$ and identity (ii) of Lemma 3.12 in \cite{Monteb}.
In the case of  $\mathbb{H}^n$ this identity
can be rewritten as follows:
\begin{equation*}\label{orie}\left\langle\gc_{\TB_{2n+1}}{\nn},\tau_j\right\rangle=\tau_j(\varpi)+\frac{1}{2}\left\langle
C^{2n+1}\cc{\nn},\tau_j\right\rangle -\left \langle
C\TB_{2n+1},\tau_j\right\rangle\qquad\mbox{for every}\quad j\in
I\ss,\end{equation*}where
$\TB_{2n+1}:=\tau_{2n+1}-\varpi\tau_1=T-\varpi\nn$. Therefore
\begin{equation*}\left\langle\gc_{T}{\nn},\tau_j\right\rangle
 -\varpi\left\langle\gc_{\nn}{\nn},\tau_j\right\rangle=\tau_j(\varpi)+\frac{1}{2}\left\langle
C^{2n+1}\cc{\nn},\tau_j\right\rangle +\varpi^2\left\langle
C^{2n+1}\cc\nn,\tau_j\right\rangle\qquad \mbox{for every}\quad
j\in I\ss.\end{equation*} By using (iii) of Lemma 3.13 in
\cite{Monteb}, we also get that
$\left\langle\gc_{T}{\nn},\tau_j\right\rangle=\frac{1}{2}\left\langle
C^{2n+1}\cc{\nn},\tau_j\right\rangle $ and the thesis easily
follows.
\end{proof}

Note that, as a byproduct of identity \eqref{idder} and  Lemma \ref{lppl}, we immediately
get the next:\begin{lemma}\label{tyie}Let $S\subset\mathbb{H}^n$ be a $\cont^2$-smooth
non-characteristic hypersurface. Then $$\left\langle{\rm Hess}\cc \phi\,
\,\nn,\nn\right\rangle=\frac{\partial^2\phi}{\partial\nn^2}+\left
\langle\frac{\qq\varpi}{\varpi},\qq
\phi\right\rangle-\varpi\frac{\partial\phi}{\partial\nn^\perp}$$for every
 $\phi\in\cin(\mathbb H^n)$.\end{lemma}

\subsection{Horizontal integration by parts}\label{IBPAA}
First, let  $S\subset\mathbb{H}^n$ be a $\cont^2$-smooth
non-characteristic hypersurface.
\begin{Defi}[Horizontal tangential operators]\label{Deflh}Let $\lg:\XX^1\ss\longrightarrow\cont(S)$  be
the 1st order
differential operator given by
\begin{eqnarray*}\lg(X):=\div\ss X + \varpi\langle C^{2n+1}\cc\nn,
X\rangle=\div\ss X -\varpi\langle\nn^{\perp}, X\rangle\qquad
\mbox{for every}\,\,X\in\XX^1\ss.\end{eqnarray*}Furthermore, let $\lh:\cont\ss^2(S)\longrightarrow\cont(S)$ be the 2nd order differential operator defined as
\begin{eqnarray*}\lh\varphi:=\Delta\ss\varphi -\varpi\frac{\partial\varphi}{\partial\nn^\perp}
\qquad\mbox{for every}\,\,\varphi\in\cont^2\ss(S).\end{eqnarray*}
\end{Defi}

Note that $\lg(\varphi X)=\varphi\lg X +
\langle\qq \varphi, X\rangle$   for every $X\in\XX^1\ss$ and every $\varphi\in\cont\ss^1(S)$. Moreover  $\lh \varphi=\lg(\qq
\varphi)$ for every $\varphi\in\cont\ss^2(S)$.

The previous definition is motivated by
Theorem 3.17 in \cite{Monteb}. For the sake of simplicity, let us
illustrate  the case of the 1st Heisenberg group
$\mathbb{H}^1$.
\begin{oss} Let $S\subset \mathbb{H}^1$ be a smooth non-characteristic surface. The $\HH$-perimeter form $\perhu$ is given by $\perhu\res
S=(\nn^\perp)^\ast\wedge\theta|_S$, where
$(\nn^\perp)^\ast:=-{\nn}_{_{\!Y}}dx + {\nn}_{_{\!X}} dy$. Note that since $\HS$ is 1-dimensional,  $\nn^\perp$ turns out to be the unique horizontal tangent
direction along $S$. Now let us
compute the exterior derivative of the 1-form $(X\LL\perhu)|_S$, for any $X\in\XX^1\ss$. Assuming
$X=f\nn^\perp$, for some smooth  $f:S\longrightarrow\R$,
yields\[d(X\LL\perhu)|_S=d\left(X\LL\left((\nn^\perp)^\ast\wedge\theta\right)\right)\Big|_S=d(f\theta)|_S=\left(df\wedge\theta
+ fd\theta\right)|_S.\]Since $d\theta=-dx\wedge dy$, by means of a
linear change of variables, we  immediately get that
 \[d\theta|_S=-(\nn)^\ast\wedge
(\nn^\perp)^\ast\big|_S=-\varpi\perhu|_S.\]Finally, since
$(\nn)^\ast\wedge\theta|_S=0$, we obtain
\[d(X\LL\perhu)|_S=\left(\nn^\perp f - f \varpi\right)\perhu|_S=\lg(X)\,\perhu|_S.\]
\end{oss}
Now let $\partial S$ be a (piecewise) $(2n-1)$-dimensional
$\cont^1$-smooth manifold, oriented by its unit normal vector
$\eta\in\TT S$ and let $\sigma^{2n-1}\rr$ be the Riemannian measure on $\partial S$,
  i.e.  $\sigma^{2n-1}\rr\res{\partial S}=(\eta\LL\sigma^{2n}\rr)|_{\partial S}$. If $X\in\cont(S, \TS)$, then   $(X\LL\perh)|_{\partial S}=\langle X, \eta\rangle
  |\PH\nu|\, \sigma^{2n-1}\rr\res{\partial S}$. The
  {\it
characteristic set} $C_{\partial S}$ of ${\partial S}$ is defined as $ C_{\partial
S}:=\{p\in{\partial S}: |\PH\nu||\P\ss\eta|=0\}$. The {\it unit
$\HS$-normal} along $\partial S$ is given by
$\eta\ss:=\frac{\P\ss\eta}{|\P\ss\eta|}$ and we may
define a homogeneous measure ${\sigma^{2n-1}\ss}$ along $\partial
S$  by setting
 ${{\sigma^{2n-1}\ss}}\res{\partial S}:=
 \left(\eta\ss\LL\perh\right)|_{\partial S}$. As for the $\HH$-perimeter, the measure ${\sigma^{2n-1}\ss}$ can be
represented by means of the Riemannian measure ${\sigma^{2n-1}\rr}$ and
it turns out that $\sigma^{2n-1}\ss\res{\partial S}=
|\PH\nu|\,|\P\ss\eta|\,\sigma^{2n-1}\rr\res{\partial S}$.\\

The following horizontal integration by parts formula can be found in
\cite{Monteb}.
\begin{teo}\label{GD}Let $S\subset\mathbb{H}^n$ be a $\cont^2$-smooth
compact non-characteristic hypersurface with piecewise
$\cont^1$-smooth boundary $\partial S$. Then
\[\int_{S}\lg(X)\,\perh=-\int_{S}\MS\langle X, \nn\rangle\,\perh +
\int_{\partial S}\langle
X,\eta\ss\rangle\,{{\sigma^{2n-1}\ss}}\qquad\mbox{for
every}\,\,X\in\XX^1\cc.\]\end{teo}

\begin{oss}\label{cpoints}
Let $S\subset \mathbb{H}^n$ be $\cont^2$-smooth compact
hypersurface with piecewise $\cont^1$-smooth boundary $\partial
S$. As already said, in this case it turns out that $\dim_{\rm
Eu-Hau}({\rm Car}_S)\leq n$. Just for the case $n=1$, we will further
suppose that {\rm $C_S$ is contained in a finite union of
$\cont^1$-smooth horizontal curves}. Under these assumptions,  one
can show that there exists a family $\{\UU_\epsilon\}_{\epsilon\geq 0}$
of open subsets of $S$ with piecewise $\cont^1$-smooth boundaries
such that:\begin{itemize}\item[{\rm(i)}] ${\rm
Car}_S\Subset\UU_\epsilon$ for every $\epsilon>0$;
\item[{\rm(ii)}] $\sigma^{2n}\rr(\UU_\epsilon)\longrightarrow 0$
for $\epsilon\rightarrow
0^+$;\item[{\rm(iii)}]$\int_{\partial\UU_\epsilon}|\PH\nu|\,\sigma^{2n-1}\rr\longrightarrow
0$ for $\epsilon\rightarrow 0^+$.\end{itemize}By using {(iii)} above,
we infer that
$\sigma^{2n-1}\ss(\partial\UU_\epsilon)\longrightarrow 0$ for
$\epsilon\rightarrow 0^+$. Therefore, in order to extend Theorem \ref{GD} to
the case $C_S\neq\emptyset$, it is then sufficient to apply it to
the non-characteristic hypersurface $S_\epsilon:=S\setminus
\UU_\epsilon$. We have
\begin{eqnarray}\label{divcar} \int_{S_\epsilon}\lg(X)\,\perh =
-\int_{S_\epsilon}\MS\langle X, \nn\rangle\,\perh +\int_{\partial
S}\langle
X,\eta\ss\rangle\,{\sigma^{2n-1}\ss}-\int_{\partial\UU_\epsilon}\langle
X,\eta\ss\rangle\,{\sigma^{2n-1}\ss}.\end{eqnarray}Note that  $\div\ss
X=\div\cc X-\langle(\mathcal{J}\cc X)\nn, \nn\rangle\in\cont(S\setminus C_S)\cap L^\infty(S)$. Furthermore, since $\varpi\in L^1(S, \perh)$, by using
 dominated convergence and {(ii)}, we get that
$$\lim_{\epsilon\rightarrow
0^+}\int_{S_\epsilon}\left(\div\ss X + \langle C\cc(\varpi)\nn,
X\rangle\right)\perh=\int_S\left( \div\ss X + \left\langle
C\cc(\varpi)\nn, X\right\rangle\right)\perh.$$It is
not difficult\footnote{Indeed, note that $|\MS|=\frac{\left|\div\cc(\PH\nu)-\langle\grad\cc|\PH\nu|,\nn\rangle\right|}{|\PH\nu|}$.} to show that $\MS\in L^1(S, \perh)$.  Therefore\[\int_{S}\MS\langle X,
\nn\rangle\,\perh=\lim_{{\epsilon\rightarrow
0^+}}\int_{S_\epsilon}\MS\langle X, \nn\rangle\,\perh.\] Finally,
by applying {(iii)}, we see that the integral over
$\partial\UU_\epsilon$ in \eqref{divcar} goes to $0$  as long as
$\epsilon\rightarrow 0^+$. This  shows the validity of Theorem \ref{GD} even if $C_S\neq\emptyset$. An analogous argument can be used, in order to prove the validity of some Green-type formulas up to $C_S$, at least if we consider functions belonging to $\cont^2\ss(S)$.
\end{oss}
We now collect some useful formula concerning the operator  $\lh$.

\begin{teo}\label{Properties}Let $S\subset\mathbb{H}^n$ be a $\cont^2$-smooth compact
hypersurface with piecewise $\cont^1$-smooth boundary $\partial
S$.  If $n=1$, assume further that $C_S$ is contained in a finite
union of $\cont^1$-smooth horizontal curves. Then, the following
hold:\begin{itemize} \item[{\rm(i)}] $\int_{S}\lh
\varphi\,\perh=0$ for every compactly supported
$\varphi\in\cont^2\ss(S)$; \item[ {\rm(ii)}]$\int_{S}\lh
\varphi\,\perh=\int_{\partial
S}{\partial\varphi}/{\partial\eta\ss}\,\sigma^{2n-1}\ss$ for every
$\varphi\in\cont^2\ss(S)$;\item[{\rm(iii)}]$\int_{S}\psi\,\lh
\varphi\,\perh=\int_{S}\varphi\,\lh\psi \,\perh$  for every
compactly supported $\varphi,\,\psi\in\cont^2\ss(S)$;
\item[{\rm(iv)}]$\int_{S}(\psi\,\lh \varphi-\varphi\,\lh\psi)
\,\perh=\int_{\partial
S}(\psi{\partial\varphi}/{\partial\eta\ss}-\varphi{\partial\psi}/{\partial\eta\ss})
\,\sigma^{2n-1}\ss$ for every
$\varphi,\,\psi\in\cont^2\ss(S)$;\item[
{\rm(v)}]$\int_{S}\psi\,\lh
\varphi\,\perh=-\int_{S}\langle\qq\varphi,\qq\psi\rangle\,\perh +
\int_{\partial
S}\psi{\partial\varphi}/{\partial\eta\ss}\,\sigma^{2n-1}\ss$
 for every
$\varphi,\,\psi\in\cont^2\ss(S)$;\item[{\rm(vi)}]$\int_{S}\lh
(\varphi^2)\,\perh=2\int_{S}\varphi\lh
\varphi\,\perh+2\int_{S}|\qq\varphi|^2\,\perh =
\int_{\partial S}{\partial\varphi^2}/{\partial\eta\ss}\,\sigma^{2n-1}\ss$
for every $\varphi\in\cont^2\ss(S)$.
\end{itemize}
\end{teo}\begin{proof}See \cite{Monteb} for the non-characteristic case. If
$C_S\neq\emptyset$, the proof follows  by arguing exactly as in Remark \ref{cpoints}.
\end{proof}

We end this section by stating a first
 consequence.

\begin{Prop}[Hopf's Lemma]\label{hopf}Let $S\subset\mathbb{H}^n$
 be a $\cont^2$-smooth closed hypersurface.
If $n=1$, assume further that $C_S$ is contained in a finite union
of $\cont^1$-smooth horizontal curves.  If $\varphi$ is a
$\cont^2\ss$-smooth function  such that
 $\lh\varphi\geq 0$ everywhere in $S\setminus C_S$, then $\varphi$ is a constant function and  $\lh\varphi=0$.\end{Prop}
\begin{proof}Since $\partial S=\emptyset$, by using (iii) in Theorem \ref{Properties}, we have
$\int_S\lh\varphi=0$ and since $\lh\varphi\geq 0$, it follows that
$\lh\varphi=0$ on $S\setminus C_S$. Thus we get that $\int_{S}|\qq\varphi|^2\perh=0$ and hence
$|\qq\varphi|=0$. Now we claim that $\varphi$ must be constant. If $n>1$, this
follows by  the {\it bracket-generating
condition} (i.e. $[\HS,\HS]=\TS$) which is satisfied by the
horizontal subbundle $\HS$ of any smooth  non-characteristic
hypersurface $S\subset\mathbb{H}^n$. Note  that
under our assumptions, $S\setminus C_S$ turns out to be a finite
union of non-characteristic open and connected subsets of $S$.
Hence, the claim
follows by the continuity of $\varphi$. If $n=1$, arguing as above implies that $\varphi$ is
constant along each leaf of the so-called {\rm characteristic
foliation} of $S$ and the thesis easily follows.
\end{proof}

\subsection{Further remarks about the horizontal integration by parts up to $C_S$}\label{Sez03}

It is well-known that Stokes formula is concerned with integrating a $k$-form over a $k$-dimensional manifold with boundary.
A common way to state this fundamental result is the following:

\begin{Prop}[Stokes formula]\label{ST} Let $M$ be an oriented $k$-dimensional manifold of class $\cont^2$ with boundary
$\partial M$. Then $$\int_M d\alpha=\int_{\partial M}\alpha$$for every compactly supported $(k-1)$-form $\alpha$ of class $\cont^1$.
\end{Prop}
One requires $M$ to be of class $\cont^2$ for a technical reason concerning \textquotedblleft pull-back\textquotedblright of differential forms.
Without much effort, it is possible to extend Proposition \ref{ST}
to the case where:\begin{center}$(\star)$\quad \it $M$ is of class $\cont^1$ and $\alpha$ is a $(k-1)$-form such that $\alpha$ and $d\alpha$ are continuous.

\end{center}
 \rm
For a more detailed discussion we refer the reader to the book by Taylor \cite{Taylor}. Nevertheless, it is worth noting that
much more general versions of Stokes formula are available in literature (see, for instance, \cite{FE}) and that researches aiming to generalize it are still intense.

In our setting, we have already discussed horizontal integration by parts formulas. However, the  previous condition $(\star)$ can easily be used in order to extend our formulas to vector fields (and functions)  possibly singular at the characteristic set $C_S$. More precisely, let $S\subset\mathbb{H}^n$
 be a $\cont^2$-smooth hypersurface  with piecewise $\cont^1$-smooth boundary $\partial S$ and let $X\in\cont^1(S\setminus C_S, \HS)$. Furthermore, set $$\alpha_X:=(X\LL \perh)|_S.$$
Then, condition $(\star)$ simply requires that $\alpha_X$ and $d\alpha_X$ be continuous on $S$. Note that $X$ is of class $\cont^1$ out of $C_S$ but may be singular at $C_S$.
\begin{Defi}Let $X\in\cont^1(S\setminus C_S, \HS)$ and set
 $\alpha_X:=(X\LL \perh)|_S$. We say that $X$ is \rm admissible (for the horizontal divergence formula)  \it if, and only if, the differential forms $\alpha_X$ and $d\alpha_X$ are continuous on all of $S$. Furthermore, we say that $\phi\in\cont^2\ss(S\setminus C_S)$  is \rm admissible \it if, and only if, $\qq \phi$ is admissible (for the horizontal divergence formula).
\end{Defi}

In other words, we have the formulas:\begin{itemize}
 \item $\int_{S}\lg(X)\,\perh=
\int_{\partial S}\langle
X,\eta\ss\rangle\,{{\sigma^{2n-1}\ss}}$ for every \it admissible \rm $X\in\cont^1(S\setminus C_S, \HS)$;
 \item $\int_{S}\lh \phi\,\perh=
\int_{\partial S}\langle
\qq\phi,\eta\ss\rangle\,{{\sigma^{2n-1}\ss}}$ for every \it admissible \rm $\phi \in\cont^2\ss(S\setminus C_S, \HS)$.
\end{itemize}

For later purposes, we have now to define a space of \textquotedblleft admissible\textquotedblright   functions for the horizontal Green's
 formulas (iii)-(vi) of Theorem \ref{Properties}.
\begin{Defi}\label{adm} We say that  $\phi\in\cont^2\ss(S\setminus C_S)$  is \rm admissible for the horizontal Green's
 formulas \it if, and only if, $\psi\,\qq \phi$ is admissible (for the horizontal divergence formula)  for every $\psi\in\cont^2\ss(S\setminus C_S)$ such that  $\psi\,\qq \psi$ is admissible (for the horizontal divergence formula). We shall denote by $\varPhi(S)$ the  space of all admissible functions for the horizontal Green's
 formulas.
\end{Defi}

Obviously, the formulas (iii)-(vi) in Theorem \ref{Properties} extend to functions belonging to the space $\varPhi(S)$ of all admissible functions for the horizontal Green's
 formulas.

 \subsection{Closed eigenvalue problem for $\lh$
}\label{auts}

 Let $S\subset\mathbb{H}^n$ be a compact closed
$\cont^2$-smooth hypersurface and denote by $L^2(S, \perh)$ the space
of all $\perh$-measurable functions $f$ on $S$ such that
$\int_S|f|^2\,\perh<+\infty$. Obviously,  $L^2(S, \perh)$ is a Hilbert space  with respect to the inner product $(\cdot,\cdot)$ defined by $(f, g):=\int_S
f\,g\,\perh$ for every $f, g\in L^2(S, \perh)$. Its associated
norm  is denoted by $\|\cdot\|$. For $\perh$-measurable horizontal tangent vector fields
we define an inner product  (with associated norm $\|\cdot\|_0$) by setting
\begin{equation}\label{l2n}(X, Y)_0:=\int_S\langle X, Y\rangle\,\perh \qquad\left(  \|X\|_0^2:=\int_S|X|^2\,\perh\right) .\end{equation}The resulting metric space, denoted by $\mathcal{L}^2(S, \perh)$,  identifies
with the space of $\perh$-measurable horizontal tangent vector fields for which the integral
\eqref{l2n} is finite. If $f\in\cont^1\ss(S)$ and
$X\in\XX^1\ss=\cont^1(S, \HS)$, then \[(\qq f, X)_0=-(f, \lg X).\]However, we would like to include in our analysis, more general functions.
\begin{Defi}Given a function $f\in L^2(S, \perh)$, we say that $Y\in
\mathcal{L}^2(S, \perh)$ is a $\HS$-{\rm weak derivative of $f$}
if, and only if, one has
$(Y, X)_0=-\int_S f\,\lg X\,\perh$ for every $X\in\XX^1\ss$.\end{Defi}

Since it is not difficult to show that there exists at
most one such $Y\in\mathcal{L}^2(S, \perh)$, we shall write
$Y=\qq f$. In the sequel, we will denote by $\mathcal{H}(S,
\perh)$ the subspace of ${L}^2(S, \perh)$ of functions
having $\HS$-weak derivatives. This space is  endowed
with the inner product $(f, g)_1:=(f, g) + (\qq f, \qq g)_0$ with
associated norm $\|f\|^2_1:=\|f\|^2+\|\qq f\|_0^2$.
\begin{oss}We stress that the space $\varPhi(S)$ of all admissible functions  for the horizontal Green's
 formulas is a subspace of ${\mathcal H}(S, \perh)$.
\end{oss}

The {\it
energy integral} in $\mathcal{H}(S, \perh)$ is the symmetric
bilinear form defined by\[\mathcal{E}(f, g):=\int_S\langle\qq f,
\qq g\rangle\,\perh\]for every $f, g\in\mathcal{H}(S, \perh)$.

Here we
are concerned  with the validity, under suitable assumptions, of the formula:
\begin{equation}\label{h1ipp}(\lh \phi, f)=-\mathcal{E}(\phi, f).\end{equation}
By the results of Section \ref{IBPAA}, this formula holds true
for every $\phi\in\cont^2\ss(S)$ and for every $f\in\cont\ss^1(S)$.
More generally,  let us assume  $\phi\in\varPhi(S)$; see Definition \ref{adm} in  Section \ref{Sez03}. In this case formula \eqref{h1ipp} defines a
linear functional $L_\phi$ on $\cont\ss^1(S)$ as a subspace of
$\mathcal{H}(S, \perh)$ satisfying: \[|L_\phi(f)|\leq\|\qq
\phi\|_0\|\qq f\|_0\leq\|\qq \phi\|_0\|f\|_1.\]In other words,
$L_\phi$ turns out to be a bounded linear functional on
$\cont\ss^1(S)\subset \mathcal{H}(S, \perh)$ with norm bounded from the
above by $\|\qq \phi\|_0$ and so it can be extended to a bounded
linear functional on all of $\mathcal{H}(S, \perh)$. Therefore,  formula \eqref{h1ipp} turns out to be valid
for every $\phi\in \varPhi(S)$ and every $f\in\mathcal H(S,\perh)$.

We are now in a position to formulate the closed eigenvalue problem for the
operator $\lh$.
\begin{claima}\label{PEH}Let $S\subset\mathbb{H}^n$ be a compact closed
 hypersurface of class $\cont^2$. The {\rm closed eigenvalue problem}
for the operator $\lh$ on $S$ is to find all real numbers
$\lambda$ for which there exist non-trivial solutions $\varphi\in\varPhi(S)$
to the following
\begin{displaymath}
{\rm (P)}\,\,\left\{%
\begin{array}{ll}
    \quad\,\lh\varphi=-\lambda\,\varphi \quad\mbox{on}\,\,S\setminus C_S\\
    \int_S\varphi\,\perh=0.\\
\end{array}%
\right..\end{displaymath}
\end{claima}

Some remarks are in order. \begin{itemize}\item \it If $\phi\in\varPhi(S)$ is an eigenfunction of Problem
\ref{PEH}, then its eigenvalue $\lambda$ must be non-negative. \rm The
proof of this claim is based on the identity
$\int_{S}\left(\phi\,\lh\phi+|\qq
\phi|^2\right)\,\perh=0,$ from which we get
\begin{eqnarray*}\int_{S}|\qq
\phi|^2\,\perh=-\int_{S}\phi\,\lh\phi\,\perh=\lambda\,\int_{S}\phi^2\,\perh.\end{eqnarray*}
Note that  if  $\lambda=0$,  it follows that $\phi$ must be
constant along $S$; see Proposition \ref{hopf}. Hence
$\lambda_0=0$ cannot be an eigenvalue of Problem \ref{PEH}.

\item \it Eigenspaces belonging to different eigenvalues are orthogonal in
$L^2(S, \perh)$. \rm This follows basically
from  formula (iii) in Theorem \ref{Properties}, i.e.
\[\int_S\left(\phi\lh \psi-\psi\lh \phi\right)\,\perh=0.\]Indeed, let $\phi,\,\psi$ be eigenfunctions of
$\lambda$ and $\tau$, respectively. Then\[\int_S\left(\phi\lh
\psi-\psi\lh \phi\right)\,\perh=(\lambda-\tau)\int_S
\phi\,\psi\,\perh=0.\]
\item The dimension of each eigenspace is called
the {\it multiplicity} of the eigenvalue. Later on, we shall list
the eigenvalues as
$0< \lambda_1\leq\lambda_2\leq\ldots\nearrow +\infty,$ with each
eigenvalue repeated according its multiplicity.

\end{itemize}

\begin{teo}[Rayleigh's Theorem]Let $S\subset \mathbb{H}^n$ be
a $\cont^2$-smooth compact closed hypersurface and let us consider the
closed eigenvalue problem for the operator $\lh$ on $S$; see Problem \ref{PEH}. Assume
that there exist eigenvalues
\begin{equation}\label{list}0<
\lambda_1\leq\lambda_2\leq\ldots\nearrow
+\infty,\end{equation}where
 each eigenvalue is repeated the number of times equal to its multiplicity. Then, for
 every non-zero function
 $f\in \mathcal{H}(S, \perh)$ one has $\lambda_1\leq \frac{\mathcal{E}(f, f)}{\|f\|^2}$, with
  equality if, and only if, $f$ is an eigenfunction of $\lambda_1$. If $\{\phi_1, \phi_2,\ldots\}$
 is a complete orthonormal basis of $L^2(S, \perh)$ such that $\phi_j$ is an eigenfunction of $\lambda_j$
 for every $j=1, 2,...$, then for any non-zero function $f\in \mathcal{H}(S, \perh)$ such that $(f, \phi_1)
 = \ldots=(f, \phi_{k-1})=0$, it turns out that $\lambda_k\leq \frac{\mathcal{E}(f, f)}{\|f\|^2}$,
 with equality if, and only if, $f$ is an eigenfunction of $\lambda_k$.
 \end{teo}

\begin{proof}The proof is based on the validity of formula \eqref{h1ipp} and it
is completely analogous to the Riemannian one, for which we refer
the reader to Chavel's book \cite{Ch3}.\end{proof}

\begin{teo}[Max-Min Theorem] Let $\eta_1,...,\eta_{k-1}\in L^2(S, \perh)$ and let $\mu=\inf\frac{\mathcal{E}(f, f)}{\|f
\|^2}$, where $f$ belongs to the subspace, less the origin, of
functions in $\mathcal{H}(S, \perh)$ orthogonal to
$\eta_1,...,\eta_{k-1}$ in $L^2(S, \perh)$. Then, for eigenvalues
given in \eqref{list}, we have $\mu\leq \lambda_k$. If
$\eta_1,...,\eta_{k-1}$ are orthonormal and each $\eta_l$ is an
eigenfunction of $\lambda_l$ for every $l=1,...,k-1$, then $\mu=\lambda_k$.
\end{teo}

\begin{proof}The proof can be done by repeating the arguments of the Riemannian case; see \cite{Ch3}.\end{proof}

 \section{Closed eigenvalue problem for $\lh$ on Isoperimetric Profiles}

\subsection{Heisenberg Isoperimetric Profiles}\label{Sez3}

This section is devoted to study some features of the unit
Isoperimetric Profile $\Sph\subset\mathbb{H}^n$. As already said,
Isoperimetric Profiles are compact hypersurfaces, similar to
ellipsoids,  which turn out to be fibred by CC-geodesics; see also
\cite{LeoM}. Their importance comes from  a longstanding
conjecture, usually attributed to Pansu, stating that they
minimize the $\HH$-perimeter among finite $\HH$-perimeter sets
having fixed volume. In other words, the Isoperimetric Profile is
the main (perhaps only) candidate to solve the sub-Riemannian
isoperimetric problem in $\mathbb{H}^n$. There is a wide
literature on this subject and many partial answers, for which we
refer the reader to \cite{CDPT}, \cite{DanGarN8}, \cite{gar},
\cite{gar2}, \cite{LeoM}, \cite{Monti, Monti2}, \cite{Monti3},
\cite{Monti4}, \cite{Ni}, \cite{P1, Pansu2}, \cite{RR}. In
particular, in \cite{RR} it is shown that the conjecture is true
for compact
$\cont^2$-smooth surfaces in $\mathbb{H}^1$.\\
 Let us preliminary state some basics about CC-geodesics in $\mathbb{H}^n$. By definition, CC-geodesics are horizontal curves
which minimize the CC-distance $\dc$ between two given points. They turn out to be solutions of the following system of O.D.E.s:
\begin{eqnarray*}\left\{\begin{array}{ll}\dot{x}=P\cc\\\dot{P}\cc=-
P_{2n+1}C^{2n+1}\cc
P\cc\\\dot{P}_{2n+1}=0,\end{array}\right.\end{eqnarray*}or
equivalently, we have to solve $\dot{P}\cc=P_{2n+1} P^\perp\cc$,
where $P\cc=(P_1,...,P_{2n})^{\rm Tr},\,|P\cc|=1$ and $P_{2n+1}$
is a constant parameter along the solution. We stress that $P=(P\cc, P_{2n+1})\in\R^{2n+1}$ can be thought of as Lagrangian multipliers and that the previous system can
be deduced by directly minimizing a constrained Lagrangian; see
\cite{Montegv} and references therein, or \cite{Montgomery}.
Unlike the Riemannian case,
CC-geodesics in $\mathbb{H}^n$ depend not only on the initial
point $x(0)$ and on the initial direction $P\cc(0)$, but also on
the parameter $P_{2n+1}$ which is a type of curvature for the solution. Indeed, if $P_{2n+1}=0$, CC-geodesics are
Euclidean horizontal lines, while  if $P_{2n+1}\neq 0$,
every CC-geodesic $x(t)$ starting from a fixed point $x_0$ is a sort of
\textquotedblleft helix\textquotedblright\footnote{If $n=1$, $x(t)$ is a circular helix with axis
parallel to the vertical direction $T$ and whose slope depends on
$P_3$. We stress that the projection of $x(t)$ onto $\R^2$ turns
out to be a circle whose radius explicitly depends on $P_3$.}. In this case, the horizontal projection of $x(t)$ onto
$\HH_0=\R^{2n}$ belongs to a sphere whose radius only depends on
$P_{2n+1}$. In particular, $x(t)$ touches the $T$-line  through  $x_0$ infinitely many times.
However,  $x(t)$
minimizes the CC-distance from $x_0$ only on the connected subset
determined by the first point $x_1$ belonging to the
$T$-line  through $x_0$. Set
$x_0={\mathcal S},\,x_1={\mathcal N}$, and call them the North and South poles. By rotating this
curve around the $T$-line  through $x_0$, we obtain a
closed convex hypersurface of constant horizontal mean curvature
 called {\it Isoperimetric Profile}. Below we will study some features of a \textquotedblleft model\textquotedblright Isoperimetric Profile having barycenter at the identity
$0\in\mathbb{H}^n$. Of course, due to left-translations and dilations, this is not restrictive at all.

Let $\rho:=\|z\|=\sqrt{\sum_{i=1}^n (x_i^2+y_i^2)}$ be the (Euclidean) norm of
$z=(x_1,y_1,...,x_i,y_i,...,x_n,y_n)\in\R^{2n}$. Moreover, let
$u_0:\overline{B_1(0)}:=\{z\in\R^{2n}: 0\leq\rho\leq 1\}\longrightarrow\R$,
\begin{equation}\label{u0}u_0(\rho):=\frac{\pi}{8} +
\frac{\rho}{4}\sqrt{1-\rho^2}-\frac{\rho}{4}\arcsin\rho\quad(0\leq\rho\leq
 1);\end{equation}see also \cite{DanGarN8, gar, gar2}, \cite{LeoM}. Setting $\Sph^{\pm}:=\left\{p=\exp\left(\sum_{i=1}^n(x_iX_i+y_iY_i)+ tT\right)\in\mathbb{H}^n: t=\pm
 u_0\right\}$ we define the {\it Heisenberg unit Isoperimetric Profile} $\Sph$
to be the compact hypersurface obtained by gluing together $\Sph^{+}$
and $\Sph^{-}$, i.e.
 $\Sph=\Sph^{+}\cup\Sph^{-}$. Since
$\nabla_{\R^{2n}} u_0={u_0}'(\rho)\frac{z}{\rho}$, it follows that
the Euclidean unit normal ${\rm n}\eu$ along $\Sph^\pm$ is given by
${\rm n}^{\pm}\eu=\frac{\left(-\nabla_{\R^{2n}}
 u_0,\pm 1\right)}{\sqrt{1+\|\nabla_{\R^{2n}}
 u_0\|^2}}$, from which we compute the Riemannian unit normal $\nu$. More precisely, we have
\begin{eqnarray*}\nu^{\pm}=\frac{\left(-\nabla_{\R^{2n}}
 u_0 \pm \frac{z^{\perp}}{2},\pm 1\right)}{\sqrt{1+\|\nabla_{\R^{2n}} u_0\|^2+\frac{\rho^2}{4}}}.\end{eqnarray*}
Now since
$$\frac{\partial u_0}{\partial\rho}(\rho)=\frac{1}{4}\left(\sqrt{1-\rho^2}-\frac{\rho^2}{\sqrt{1-\rho^2}}-\frac{1}{\sqrt{1-\rho^2}}\right)
=\frac{-\rho^2}{2\sqrt{1-\rho^2}},$$we get that
\begin{eqnarray*}{\nu}^{\pm}\cc=\frac{\left(-\nabla_{\R^{2n}} u_0
\pm \frac{z^{\perp}}{2}\right)}{\sqrt{\|\nabla_{\R^{2n}}
u_0\|^2+\frac{\rho^2}{4}}}=z\pm\frac{\sqrt{1-\rho^2}}{\rho}z^\perp \end{eqnarray*}
and so
\begin{eqnarray*}({\nu}^{\pm}\cc)^\perp=\left(-z\pm\frac{\sqrt{1-\rho^2}}{\rho}z^\perp\right)^\perp
=z^\perp\mp\frac{\sqrt{1-\rho^2}}{\rho}z.\end{eqnarray*} Moreover,
since
$|\PH(\nu^{\pm})|=
\frac{\sqrt{\|\nabla_{\R^{2n}} u_0\|^2+\frac{\rho^2}{4}}}
{\sqrt{1+\|\nabla_{\R^{2n}}
u_0\|^2+\frac{\rho^2}{4}}},$ we also get that
\begin{eqnarray*}\perh\res\Sph^{\pm}=|\PH(\nu^{\pm})|\,\sigma^{2n}\rr\res\Sph^{\pm}
={\sqrt{\|\nabla_{\R^{2n}} u_0\|^2+\frac{\rho^2}{4}}}\,dz\res
B_1(0)=\frac{\rho}{2\sqrt{1-\rho^2}}\,dz\res
B_1(0).\end{eqnarray*}Furthermore $\MS=-\div\cc\nn=-2n$. By its own definition, we have
\begin{eqnarray*}\varpi^{\pm}=\frac{\pm 1}{\|\nabla_{\R^{2n}} u_0+ \frac{z^{\perp}}{2}\|}=\frac{\pm 1}
{\sqrt{\|\nabla_{\R^{2n}} u_0\|^2+\frac{\rho^2}{4}}}=\frac{\pm
1}{\sqrt{{\frac{\rho^4}{4(1-\rho^2)}+ \frac{\rho^2}{4}}}}=\pm
2\frac{\sqrt{1-\rho^2}}{\rho}.\end{eqnarray*}Note that  $\dg\varphi\equiv\nabla_{\R^{2n}}\varphi$
for every function $\varphi:\mathbb H^n: \longrightarrow\R$ independent of $t$. So we  get
\begin{eqnarray*}\dg\varpi^\pm=\nabla_{\R^{2n}}\varpi^\pm=
\pm\frac{\partial}{\partial\rho}\left(2\frac{\sqrt{1-\rho^2}}{\rho}\right)\frac{z}{\rho}
=\mp
\left(\frac{2}{\rho^2\sqrt{1-\rho^2}}\right)\frac{z}{\rho}\end{eqnarray*}and hence
\begin{eqnarray}\label{derombar}\frac{\partial\varpi}{\partial\nn^\perp}=
\mp\frac{2}{\rho^2\sqrt{1-\rho^2}}\left\langle\frac{z}{\rho},({\nu}^{\pm}\cc)^\perp\right\rangle=\frac{1}{\rho\sqrt{1-\rho^2}}\frac{2\sqrt{
1-\rho^2}}{\rho}=\frac{2}{\rho^2}.\end{eqnarray}

These computations will be used  throughout the next sections in order to  study the
action of the 2nd order operator $\lh$ on smooth functions defined
on the unit Isoperimetric profile $\Sph$.

\subsection{Radial case and hypergeometric solutions}\label{SEZPr}

 In this paper, we will study Problem \ref{PEH} in the case of the model
Isoperimetric Profile $\Sph$. By applying previous computations, we have
\begin{displaymath}\lh\varphi=\Delta\ss\varphi-2\frac{\sqrt{1-\rho^2}}
    {\rho}\,\left\langle\dg\varphi,\left(z^\perp-z\frac{\sqrt{1-\rho^2}}
    {\rho}\right)\right\rangle;\end{displaymath}see Section
    \ref{Sez3}. In the general case, we have to find
$\cont^2\ss$-smooth solutions  on
$\Sph\setminus\{\mathcal{N}, \mathcal{S}\}$ which belong to the space $\varPhi(\Sph)$ of all admissible functions  for the horizontal Green's
 formulas; see Section \ref{Sez03}. In particular, these functions must belong
to the horizontal tangent Sobolev space $\mathcal{H}(S, \perh)$.

Let $\varphi:\Sph\longrightarrow\R$ be the
restriction of a function
$\widetilde{\varphi}:\mathbb{H}^n\longrightarrow\R$, i.e.
$\varphi:=\widetilde{\varphi}|_{\Sph}$, and fix Euclidean
cylindrical coordinates on $\mathbb{H}^n$ , i.e. $(\rho, \xi,
t)\in\R_+ \times \mathbb{S}^{2n-1}\times \R$, where
$\rho=|z|$ and $\xi=\frac{z}{|z|}$. Hereafter, the exponential
coordinates $(z, t)\in\R^{2n+1}$  of  each point $p\in\mathbb{H}^n$ will be
identified with the triple $(\rho, \xi, t)$.

\begin{oss}Since $\Sph$ is the union
of two $T$-graphs, if $p\in\Sph$ then either $p\in\Sph^+$  or
$p\in\Sph^-$. If
 $p=(z, t)\in\Sph^\pm$, one must have $t=\pm u_0(\rho)$; see \eqref{u0}. Using spherical coordinates
 $(\rho, \xi)\in\R_+ \times \mathbb{S}^{2n-1} $
 on $\R^{2n}$, we get that any real valued function defined on  $\Sph^\pm$ can be seen as a function of
 the variables $(\rho, \xi)\in[0, 1]\times\mathbb{S}^{2n-1}$. In the sequel,
 for every $\varphi:\Sph\longrightarrow\R$
 we shall set $\varphi^\pm:=\varphi|_{\Sph^\pm}$
 and assume that the restrictions $\varphi^\pm:[0,
 1]\times\mathbb{S}^{2n-1}\longrightarrow\R$
satisfy the ``compatibility'', or continuity, constraint: $\varphi^+(1, \xi)=\varphi^-(1, \xi)$ for every $\xi\in \mathbb{S}^{2n-1}$.
\end{oss}
Because the radial symmetry of
$S=\Sph$, in this section we shall preliminarily study the case of radial functions.
Below, radial derivatives will be denoted by
$\varphi':=\frac{\partial\varphi}{\partial \rho}$ and
$\varphi'':=\frac{\partial^2\varphi}{\partial \rho^2}$.

\begin{no}
We set $g\cc:=\langle z,\nn\rangle$ and $g\cc^\perp:=\langle
z,\nn^{\perp}\rangle$. The function $g\cc$ is called {\rm horizontal
support function} associated with $\Sph$.\end{no}

 Under the
previous assumptions, we first compute
\begin{eqnarray*}\Delta\ss\varphi&=&\sum_{i\in I\ss}\tau_i(\tau_i(\varphi))=
\sum_{i\in I\ss}\tau_i\langle\dg\varphi,\tau_i\rangle=\sum_{i\in
I\ss}\tau_i\left\langle\varphi'\frac{z}{\rho},\tau_i\right\rangle\\&=&\sum_{i\in
I\ss}\left(\left\langle\gc_{\tau_i}\left(\varphi'\frac{z}{\rho}\right),
\tau_i\right\rangle+\left\langle\varphi'\frac{z}{\rho},\gc_{\tau_i}\tau_i\right\rangle\right)
\\&=&\MS \left(\frac{\varphi'}{\rho}\right)g\cc + \sum_{i\in I\ss}\tau_i\left(\frac{\varphi'}{\rho}\right)\langle
z,\tau_i\rangle+  \left(\frac{\varphi'}{\rho}\right)\sum_{i\in
I\ss}\langle\gc_{\tau_i}z,\tau_i\rangle\\&=&\MS
\left(\frac{\varphi'}{\rho}\right)g\cc  + \sum_{i\in
I\ss}\left\langle\dg\left(\frac{\varphi'}{\rho}\right),\tau_i\right\rangle\langle
z,\tau_i\rangle +
(\underbrace{2n-1}_{=\dim\HS})\left(\frac{\varphi'}{\rho}\right)\\&=&
\left(\frac{\varphi'}{\rho}\right)(\MS g\cc +  2n-1) +
\frac{\varphi''\rho-\varphi'}{\rho^2}\sum_{i\in
I\ss}\left\langle\frac{z}{\rho},\tau_i\right\rangle\langle
z,\tau_i\rangle\\&=& \left(\frac{\varphi'}{\rho}\right)(\MS g\cc +
2n-1) + \frac{\varphi''\rho-\varphi'}{\rho^3}\sum_{i\in
I\ss}\langle z,\tau_i\rangle^2\\&=&
\left(\frac{\varphi'}{\rho}\right)(\MS g\cc + 2n-1) +
\frac{\varphi''\rho-\varphi'}{\rho^3}(\rho^2-g\cc^2).
\end{eqnarray*}Using the last computation together with the identity $\dg\varphi=\varphi'\frac{z}{\rho}$
yields\begin{eqnarray*}\lh\varphi&=&\left(\frac{\varphi'}{\rho}\right)(\MS
g\cc + 2n-1) +
\frac{\varphi''\rho-\varphi'}{\rho^3}(\rho^2-g\cc^2)+2\varphi'\frac{{1-\rho^2}}
    {\rho}\\&=&
\frac{(\varphi''\rho-\varphi')}{\rho}\left(1-\left(\frac{g\cc}{\rho}\right)^2\right)+\frac{{\varphi'}}{\rho}(\MS
g\cc + 2n+1-2\rho^2).
   \end{eqnarray*}Since  $g\cc=\rho^2$ and $\MS=-2n$,
we get that
\[\lh\varphi=
\left(\varphi''-\frac{\varphi'}{\rho}\right)(1-{\rho}^2)+\frac{\varphi'}{\rho}\left((Q-1)
-{Q}{\rho^2}\right)\\\label{carra}=
\varphi''(1-{\rho}^2)+\frac{\varphi'}{\rho}\left(2n
-(2n+1){\rho^2}\right).
   \]

The last one is an interesting ordinary 2nd order differential operator
of hypergeometric type. Hence, under our current assumptions, Problem \ref{PEH} reduces to study a 1-dimensional eigenvalue problem\begin{eqnarray}\label{azcaazca1}
\varphi''(1-{\rho}^2)+\frac{\varphi'}{\rho}\left(2n
-(2n+1){\rho^2}\right)=-\lambda\varphi\qquad
(\lambda\in\R_+ )
   \end{eqnarray} subject to the integral condition
$\int_{\Sph}\varphi\,\perh=0$.\\

\begin{oss}Instead
of $\varphi$, we will study its restrictions $\varphi^\pm$  to  the
North and South hemispheres $\Sph^\pm$ and suppose that $\varphi^\pm$ are  smooth
  solutions  to \eqref{azcaazca1}  on  $]0, 1]\subset\R$, i.e. $\varphi^\pm\in\cont^2(]0, 1])$.
Moreover, $\varphi^\pm$  must satisfy the continuity constraint, i.e. $\varphi^+(1, \xi)=\varphi^-(1, \xi)$ {for every}
$\xi\in \mathbb{S}^{2n-1}$, together with the integral condition
$\int_{\Sph}\varphi\,\perh=0$. The last one can be
translated into an equation involving 1-dimensional weighted
integrals of $\varphi^\pm$.  More precisely, by using the symmetry
of $\Sph$ with respect to the horizontal hyperplane $t=0$, we get
that
\[\int_{\Sph} \varphi\,\perh=\int_{\Sph^+}
\varphi^+\,\perh+\int_{\Sph^-} \varphi^-\,\perh.\] Since
\[\perh\res\Sph^{\pm}=\frac{\rho}{2\sqrt{1-\rho^2}}\,dz\res B_1(0),\] we
obtain
\[\int_{\Sph}
\varphi\,\perh=\int_{\Sph^+} \varphi^+\,\perh+\int_{\Sph^-}
\varphi^-\,\perh=\frac{O_{2n-1}}{2}\int_0^1 \left( \varphi^++\varphi^-\right) (\rho)\frac{\rho^{2n}}{\sqrt{1-\rho^2}}\,d\rho=0.\]
It follows that any radial function
$\varphi:\Sph\longrightarrow\R$ belongs to $\mathcal{H}(S, \per)$
if, and only if, its restrictions $\varphi^\pm$ satisfy:
\begin{eqnarray}\label{Hsobolev}\varphi^\pm\in L^2\left([0, 1], \frac{\rho^{2n}}{2\sqrt{1-\rho^2}}\,d\rho
\right),\qquad (\varphi^\pm)'\in L^2\left([0, 1],
\frac{\rho^{2n}}{2\sqrt{1-\rho^2}}\,d\rho
\right).\end{eqnarray}\end{oss}

Hence, in the radial case, we can reformulate our closed
eigenvalue problem as follows:
\begin{claima}[Radial version]\label{prad}
Find all positive real numbers $\lambda\in\R_+ $ such that
there exist functions $\varphi^+, \varphi^-\in\cont^2(]0, 1])$
satisfying \eqref{Hsobolev}, which are non trivial solutions
to:\begin{displaymath}
\left\{%
\begin{array}{ll}
    {\varphi}''(1-{\rho}^2)+\frac{{\varphi}'}{\rho}\left(2n
-(2n+1){\rho^2}\right)=-\lambda{\varphi}\qquad
(\lambda\in\R_+ )\\\varphi^+(1)=\varphi^-(1)\\
    \int_0^1(\varphi^++\varphi^-)(\rho)\frac{\rho^{2n}}{\sqrt{1-\rho^2}}\,d\rho=0.
\end{array}%
\right.\end{displaymath}
\end{claima}

This one can be regarded as a Sturm-Liouville
problem with mixed conditions for a hypergeometric O.D.E. The first step is finding the general integral of equation
\eqref{azcaazca1}.

\begin{no}\label{y}Later on we will set
$${\mathbf F}(a,b,c,z):=\sum_{k=0}^{+\infty}\frac{(a)_k(b)_k}{k!(c)_k}z^k \qquad\left(
z\in\mathbb{C},\,\, |z|<1\right) $$where
$(d)_k:=\frac{\Gamma(d+k)}{\Gamma(d)}=d(d+1)(d+2)...(d+k-1)$ and
$\Gamma$ denotes the {\it Euler Gamma} function. The series ${\mathbf F}(a,b,c,z)$ is the so-called
{\rm hypergeometric series}; see \cite{GWo}.\end{no}

\begin{lemma}\label{cop}The general integral $\varphi(\rho):=\varphi(c_1,c_2, \lambda, \rho)$ of
\eqref{azcaazca1} is given by
\begin{eqnarray*}\varphi(\rho)&=&c_1{\mathbf F}\left(\frac{n-\sqrt{n^2+\lambda}}{2};
\frac{n+\sqrt{n^2+\lambda}}{2} ;\frac{1}{2}+n; \rho^2\right)\\&+&
 c_2\frac{1}{\rho^{2n-1}}{\mathbf F}\left(\frac{1-n-\sqrt{n^2+\lambda}}{2};
\frac{1-n+\sqrt{n^2+\lambda}}{2} ;\frac{3}{2}-n;
\rho^2\right).\end{eqnarray*}\end{lemma}Note that the second function is
singular at $\rho=0$. Further, set\begin{eqnarray*}\psi_1(\lambda,
\rho)&:=&{\mathbf F}\left(\frac{n-\sqrt{n^2+\lambda}}{2};
\frac{n+\sqrt{n^2+\lambda}}{2} ;\frac{1}{2}+n;
\rho^2\right),\\\psi_2(\lambda,
\rho)&:=&\frac{1}{\rho^{2n-1}}{\mathbf
F}\left(\frac{1-n-\sqrt{n^2+\lambda}}{2};
\frac{1-n+\sqrt{n^2+\lambda}}{2} ;\frac{3}{2}-n; \rho^2\right).
\end{eqnarray*}

A rigorous (and elementary) proof of the previous lemma uses
{\it Frobenius' method}. Nevertheless,  a  more concise proof can be done by using the
theory of {\it hypergeometric functions}.

\begin{proof}[Proof of Lemma \ref{cop}]First,
set $s:=\rho^2$ on $[0, 1]$ and
$\varphi(\rho):=\psi(\rho^2)=\psi(s)$ for $s\in[0, 1]$. Under this
transformation, equation \eqref{azcaazca1} becomes
\begin{eqnarray}\label{azcaazca12}
s(1-s)\psi''+{\psi'}\left(\frac{(Q-1)-Q
s}{2}\right)=-\frac{\lambda}{4}\psi.
   \end{eqnarray}Indeed, from the elementary identities $\varphi'(\rho)=2\rho\psi'(\rho^2)$,
    $\varphi''(\rho)=4\rho^2\psi''(\rho^2)+2\psi'(\rho^2)$, we
    get \[4s(1-s)\psi''+
    \left((2n+2)(1-s)-s\right)\psi'=-\lambda\psi \] and
    it is enough to divide by $4$ and to use $Q=2n+2$. In this form
    our equation can be handled by means of the general theory of
    hypergeometric functions, for which we refer the reader to the
    book by Wang and Guo \cite{GWo}. In particular, by using the
    notation of Chapter 4 in \cite{GWo}, the general
    hypergeometric equation (of the complex variable $z$) is given by
   \begin{eqnarray}\label{iper}z(1-z)\psi''+ (\gamma-(\alpha+\beta+1))\psi'-\alpha \beta \psi=0.\end{eqnarray}
This equation has two linearly independent
solutions, which can be written out in terms of hypergeometric series (see
Notation \ref{y}), i.e.

\begin{eqnarray*}\psi_1(z):={\mathbf F}(\alpha; \beta; \gamma;
z),\qquad\psi_2(z):=\frac{1}{z^{\gamma-1}}{\mathbf
F}(\alpha-\gamma+1; \beta-\gamma+1; 2-\gamma; z).
\end{eqnarray*}By analyzing the coefficients of equation
   \eqref{azcaazca12},
   we find that $\gamma=\frac{Q-1}{2}$ and so
 \begin{eqnarray*}\alpha+\beta+1=\frac{Q}{2}\qquad \alpha \beta=-\frac{\lambda}{4}.
\end{eqnarray*}Hence \[\alpha =\frac{n-\sqrt{n^2+\lambda}}{2},
\qquad \beta =\frac{n+\sqrt{n^2+\lambda}}{2},\qquad
\gamma=n+\frac{1}{2},\]and the thesis easily follows.
\end{proof}

We are now in a position to solve Problem \ref{prad}.
Remind that we are looking for solutions
$\varphi^\pm\in\cont^2(]0, 1])$ satisfying \eqref{Hsobolev}. So
let us consider the second solution $\psi_2$. It is elementary to
see that $\psi_2(\lambda, \rho)\backsim\frac{1}{\rho^{2n-1}}$ as
long as $\rho\rightarrow 0^+$ because the hypergeometric series, evaluated at
$\rho=0$, takes the value $1$. Therefore $\psi_2$ cannot satisfy the conditions in
\eqref{Hsobolev} and we get\[\varphi^\pm(\rho)=c^\pm\psi_1(\lambda, \rho)\qquad
\rho\in[0, 1],\]for some real constant $c^\pm$. We have now to use
the ``mixed'' conditions of Problem \ref{prad}. We see that the
third condition, together with the continuity of $\varphi$,
implies the following two possibilities:

\begin{itemize}\item[{\rm(i)}] $\varphi^+(\rho)=\varphi^-(\rho)$
for every $\rho\in[0, 1]$;\item[{\rm(ii)}]
$\varphi^+(\rho)=-\varphi^-(\rho)$ for every $\rho\in[0, 1]$ and
$\varphi^+(1)=\varphi^-(1)=0$.\end{itemize}

\begin{lemma}\label{ZX0}In case {(i)} holds, Problem \ref{prad} admits a countable family of eigenvalues which are given by $\lambda_{2m}:=2m(2m+
 2n),\,\,m\in\mathbb{N}$. The eigenfunction $\varphi_{2m}$ relative to $\lambda_{2m}$ is polynomial
and explicitly given, for some real constant $c$, by
\[\varphi_{2m}(\rho):=c \,{\mathbf F}\left(-m; n+m; \frac{1}{2}+n; \rho^2\right).\]

 \end{lemma}

\begin{proof}Since $\varphi^+=\varphi^-$ on $[0, 1]$, it is enough to write down the integral condition. Hence we have to find all positive real numbers
$\lambda\in\R_+ $ such that any non-singular solution $\varphi$ to
\eqref{azcaazca1} (i.e. $\varphi=c\psi_1(\lambda, \rho)$  for some
$c\in\R$) satisfies
\[
\int_0^1\varphi(\rho)\frac{\rho^{2n}}{\sqrt{1-\rho^2}}\,d\rho=0.\]By
using some classical results about hypergeometric series (or {\sc
Mathematica}) we see
 that\[\int_0^1\psi_1(\lambda,
 \rho)\frac{\rho^{2n}}{\sqrt{1-\rho^2}}\,d\rho=\frac{\sqrt{\pi}\,
 \Gamma(n+\frac{1}{2})}{2\,\Gamma(\frac{2+n-
 \sqrt{n^2+\lambda}}{2})\, \Gamma(\frac{2+n+
 \sqrt{n^2+\lambda}}{2})}.\]It remains to solve (in the variable $\lambda\geq
 0$) the following equation:\begin{equation}\label{eqr}\frac{\sqrt{\pi}\,
 \Gamma(n+\frac{1}{2})}{2\,\Gamma(\frac{2+n-
 \sqrt{n^2+\lambda}}{2})\, \Gamma(\frac{2+n+
 \sqrt{n^2+\lambda}}{2})}=0.\end{equation}From standard properties of the
Gamma function we easily get that $\frac{2+n-
 \sqrt{n^2+\lambda}}{2}$ must be a negative integer or 0, i.e.
\[\frac{2+n-
 \sqrt{n^2+\lambda}}{2}=-m,\qquad m\in\mathbb{N}.\]Therefore  $\lambda=2(m+1)(2(m+1)+ 2n)$ for any $ m\in\mathbb{N}$.
 From now on, we shall set
 \[\lambda_{2m}:=2m(2m+ 2n),\qquad m\in\mathbb{N}_+.\]
 \begin{oss}In this list, we have not included the eigenvalue $\lambda_0=0$  because, as already said, the only eigenfunction
  relative to $\lambda_0$ is $\varphi_0=0$.\end{oss}
 The eigenfunction $\varphi_{2m}$ relative to the eigenvalue
 $\lambda_{2m}\,(m\geq 1)$ is given, up to a real constant $c$, by\[\varphi_{2m}(\rho):=c\psi_1(\lambda_{2m}, \rho)=c\,
 {\mathbf F}\left(-m; n+m; \frac{1}{2}+n; \rho^2\right).\]From the general
theory of  hypergeometric series it follows that $\varphi_{2m}$ is a polynomial
function. By construction,
$\varphi_{2m}=\varphi_{2m}^+=\varphi_{2m}^-$ and we have
\[\int_0^1\varphi_{2m}(\rho)\frac{\rho^{2n}}{\sqrt{1-\rho^2}}\,d\rho=0.\]
\end{proof}

\begin{lemma}\label{ZX}In case (ii) holds,
  Problem \ref{prad} admits a countable family of eigenvalues which are given by $\lambda_{2m+1}:=(2m+1)\left((2m+1) +
 2n\right),\,\,m\in\mathbb{N}$. The eigenfunction $\varphi_{2m+1}$ relative to $\lambda_{2m+1}$ is the
 hypergeometric
function given by
\[\varphi_{2m+1}(\rho):=c \,{\mathbf F}\left(- m-\frac{1}{2}; n+m+\frac{1}{2}; \frac{1}{2}+n;
\rho^2\right),\]for some real constant $c$. Furthermore, we have
 $\varphi_{2m+1}=\varphi_{2m+1}^+$ on
${\Sph^+}$, $\varphi_{2m+1}=\varphi_{2m+1}^-$ on ${\Sph^-}$ and $\varphi_{2m+1}^+=-\varphi_{2m+1}^-$.
 \end{lemma}
\begin{proof}We are assuming that $\varphi^+=-\varphi^-$ on $[0, 1]$  and that
 $\varphi^+(1)=-\varphi^-(1)=0$. So we have to find all
positive real numbers $\lambda\in\R_+ $ such that any pair
of non-singular solutions $\varphi^\pm$ to \eqref{azcaazca1} (i.e.
$\varphi^\pm=c^{\pm}\psi_1(\lambda, \rho)$, for some $c^\pm\in\R$)
satisfies:
\[\varphi^+(1)=-\varphi^-(1)=0.\]In particular, we get that
$c^+=-c^-$ and it remains to solve (in the variable $\lambda\geq
0$) the equation $\psi_1(\lambda, \rho)=0$, i.e.
\[{\mathbf F}\Big(\frac{n-\sqrt{n^2+\lambda}}{2};
\frac{n+\sqrt{n^2+\lambda}}{2} ;\frac{1}{2}+n; 1\Big)=0.\]We
recall that a classical result about hypergeometric series states
that \[{\mathbf F}\Big(\alpha, \beta, \gamma,
1\Big)=\frac{\Gamma(\gamma)\,\Gamma(\gamma-\alpha-\beta)}{\Gamma(\gamma-\alpha)\,\Gamma(\gamma-\beta)}\]as
long as $\Re(\gamma-\alpha-\beta)>0$; see \cite{GWo} Chapter 4.7,
p. 156. Note that $\gamma-\alpha-\beta=\frac{1}{2}>0$. So we have to find the values
$\lambda$ such that
\begin{equation}\label{eqs}\frac{\Gamma(\gamma)\,\Gamma(\gamma-\alpha-\beta)}{\Gamma(\gamma-
\alpha)\,\Gamma(\gamma-\beta)}=0.\end{equation}The
general theory of the Gamma function (see \cite{GWo}, Chapter
3) says that  this equation can be solved if, and only if,
$\gamma-\alpha$ or $\gamma-\beta$ are negative integers, or 0.
This in turn implies\[\gamma-\beta=\frac{n+1-\sqrt{n^2+\lambda}}{2}=-m,\qquad
m\in\mathbb{N}.\]So we get
$\lambda_{2m+1}:=(2m+1)(2m+1+2n)$, $m\in\mathbb{N}$.
The eigenfunction $\varphi_{2m+1}$ relative to the
eigenvalue
 $\lambda_{2m+1}\,(m\in\mathbb{N})$ is given, up to a real constant
  $c$, by\[\varphi_{2m+1}(\rho):=c\psi_1(\lambda_{2m+1}, \rho)=c\,
 {\mathbf F}\left(-m-\frac{1}{2}; n+m+\frac{1}{2}; \frac{1}{2}+n;
\rho^2\right).\]The last one is a hypergeometric function and, by
construction,
  $\varphi_{2m+1}^+=-\varphi_{2m+1}^-$. Therefore
\begin{equation}\label{ovvia}\int_0^1 \left( \varphi_{2m+1}^++\varphi_{2m+1}^- \right)(\rho)
\frac{\rho^{2n}}{\sqrt{1-\rho^2}}\,d\rho=0.\end{equation} \end{proof}
 
\begin{teo}\label{RPF}Problem \ref{prad} admits a countable family of positive real eigenvalues \[\lambda_k:=k(k+2n),
\qquad k\in\mathbb{N}.\] Each eigenvalue $\lambda_k$ has
multiplicity  1 and its associated eigenfunction
$\varphi_k$, up to real constants,  is a hypergeometric function. More
precisely, we have the following:
\begin{itemize}\item[{\rm (i)}]if $k=2m$, the eigenfunction $\varphi_{2m}$
 associated with $\lambda_{2m}$ is the hypergeometric polynomial
 given by
\[\varphi_{2m}(\rho):=c \,{\mathbf F}\left(-m; n+m; \frac{1}{2}+n; \rho^2\right),\]for some real constant $c$. Moreover,
$\varphi_{2m}=\varphi_{2m}^+=\varphi_{2m}^-$ and
\[\int_0^1\varphi_{2m}(\rho)\frac{\rho^{2n}}{\sqrt{1-\rho^2}}\,d\rho=0.\]
 \item[{\rm (ii)}]if $k=2m+1$, the eigenfunction $\varphi_{2m+1}$ associated with
 $\lambda_{2m+1}$ is the hypergeometric function
\[\varphi_{2m+1}(\rho):=c\,
 {\mathbf F}\left(-m-\frac{1}{2}; n+m+\frac{1}{2}; \frac{1}{2}+n;
\rho^2\right),\]for some real
 constant $c$. Moreover, $\varphi_{2m+1}^+=-\varphi_{2m+1}^-$ and
\eqref{ovvia} holds true.\end{itemize}
\end{teo}
\begin{proof}Applying Lemma \ref{ZX0} and Lemma \ref{ZX}.\end{proof}

\begin{corollario}\label{cor1aut}Up to real constants,
 the eigenfunction $\varphi_1$ associated with the first eigenvalue $\lambda_1=Q-1$ of Problem
\ref{prad} is given by
\[\varphi_1(\rho)=\sqrt{1-\rho^2}\qquad\mbox{for every}\,\, \rho\in[0, 1].\]\end{corollario}
\begin{proof}This follows  from (ii), by using the identity
\[(1+z)^\alpha={\mathbf F}(-\alpha, \beta, \beta, -z),\]and by setting $\alpha=\frac{1}{2}$ and $z=-\rho^2$;
see, for instance, \cite{GWo}, p. 137, formula (10).\end{proof}

\begin{oss}Let $S\subset\mathbb{H}^n$ be a smooth compact
hypersurface without boundary. We stress that
\[\varphi_2(p):=(2n-1) +\MS g\cc-\varpi g\cc^\perp,\] is a zero-mean function, i.e.
 $\int_{S} \varphi_2\,\perh=0$, where
$g\cc=\langle z, \nn\rangle$ and $g\cc^\perp=\langle
z,\nn^{\perp}\rangle$. Furthermore if  $S=\Sph$, then $\varphi_2(\rho)=(Q-1)-Q\rho^2$ turns out to be radial. We also stress that, in this case, the function $\varphi_2$ is, up to real constants, the 2nd radial eigenfunction of Problem
\ref{prad}, corresponding to  $\lambda_2=2Q$.
\end{oss}

\subsection{Some remarks about the general case}\label{GENCASF}

\begin{war} From now on, when we refer to Problem \ref{PEH} it is understood that $S=\Sph$.

\end{war}
\begin{oss} As already said,
any function $\varphi:\Sph\setminus\{\mathcal{N},
\mathcal{S}\}\longrightarrow\R$ can be regarded as a restriction
of some function
$\widetilde{\varphi}:\mathbb{H}^n\longrightarrow\R$. Moreover,
$\Sph=\Sph^+\cup\Sph^-$ and the hemispheres $\Sph^\pm$ are the
 radial $T$-graphs given by  $t=\pm u_0(\rho)$, where $\rho=|z|$, $z\in
 \overline{B_1(0)}$, and
 $u_0$ was
 defined by formula \eqref{u0}. By introducing  polar coordinates\footnote{This means
that we are fixing ``ordinary'' spherical coordinates
$(\rho,\xi)\in\R_{+}\cup\{0\}\times\mathbb{S}^{2n-1}$ on
$\R^{2n}\cong\HH_0$ and so $z\in\R^{2n}$ corresponds to $(\rho,
\xi)\in\R_{+}\cup\{0\}\times\mathbb{S}^{2n-1}$.}
$(\rho,\xi,t)\in\R_{+}\cup\{0\}\times\mathbb{S}^{2n-1}\times\R$ on
$\mathbb{H}^n$ and by setting
\[\Phi^{\pm}:\overline{B_1(0)}\longrightarrow\mathbb{H}^n,\qquad
\Phi^{\pm}(\rho,\xi):=(\rho, \xi, \pm u_0(\rho)),\]it follows that
$\varphi^\pm:=\widetilde{\varphi}\circ\Phi^{\pm}$ only depends on
the variables $(\rho, \xi)$.\end{oss}
Remember also that $$\Delta\ss\varphi=\Delta\cc\varphi +
\MS\frac{\partial\varphi}{\partial\nn}-\langle{\rm
Hess}\cc\varphi\,\nn,\nn\rangle;$$see Lemma
\ref{ljjjkl}. So if
$\varphi^{\pm}=\widetilde{\varphi}\circ\Phi^{\pm}:\overline{B_1(0)}\setminus\{0\}\longrightarrow\R$,
 by a simple computation we get that:
\begin{itemize}\item $\Delta\cc\varphi^{\pm}=\Delta_{\R^{2n}}\varphi^{\pm}$,
\item $\frac{\partial^{\pm}\varphi}{\partial\nn^{\pm}}
 =\pm\langle\nabla_{\R^{2n}}\varphi^{\pm},(z\pm\kappa z^\perp)\rangle$,\item$\langle{\rm
Hess}\cc \varphi^{\pm}\,\nn,\nn\rangle=\langle{\rm Hess}_{\R^{2n}}
\varphi^{\pm}(z\pm\kappa z^\perp),(z\pm\kappa
z^\perp)\rangle$,\end{itemize}where $\nn^{\pm}=z\pm\kappa
z^\perp$ and we have set $\kappa(\rho):=\frac{\sqrt{1-\rho^2}}{\rho}$.
Because the symmetry of $\Sph$ with respect to the horizontal
hyperplane $\HH_0\cong\R^{2n}$, we will carry out the computations only
for the north hemisphere $\Sph^+\setminus\{\mathcal{N}\}$. For brevity,  set $\varphi=\varphi^+$ and $\nn=\nn^+$. We
have\begin{eqnarray*}\Delta\ss\varphi=\Delta_{\R^{2n}}\varphi -
2n\Big(\frac{\partial\varphi}{\partial
z}+\kappa\frac{\partial\varphi}{\partial z^{\perp}}\Big)
-\left\langle{\rm Hess}\cc\varphi\, \nn, \nn\right\rangle,
\end{eqnarray*}where we have used
$\MS|_{\Sph}=-2n$. Since$$\lh\varphi=\Delta\ss\varphi+2\kappa\Big(\kappa\frac{\partial\varphi}{\partial
z}-\frac{\partial\varphi}{\partial z^{\perp}}\Big),$$we obtain\begin{eqnarray}\label{1forma}\lh \varphi&=&
\Delta_{\R^{2n}}\varphi
-2(n-\kappa^2)\frac{\partial\varphi}{\partial
z}-Q\kappa\frac{\partial\varphi}{\partial
z^{\perp}}-\left\langle{\rm Hess}\cc\varphi \, \nn,
\nn\right\rangle.
\end{eqnarray}

\begin{no}Set
\[\zeta:=\frac{z^{\perp}}{\rho}\in\TT\mathbb{S}^{2n-1}.\]Hereafter, either
$\frac{\partial\varphi}{\partial\rho}$ or $\varphi'_\rho$, will
denote the 1st partial derivative of $\varphi$ with respect to the
variable $\rho=|z|$. Moreover, either $\frac{\partial\varphi}{\partial\zeta}$ or $\varphi'_\zeta$ will  denote
the 1st partial derivative of $\varphi$ with respect to the \textquotedblleft angular\textquotedblright  variable $\zeta$.
An analogous notation will be used for 2nd
partial derivatives.
\end{no}
\begin{lemma}\label{rftr}We have $\div_{\TT\mathbb{S}^{2n-1}}\zeta=0.$\end{lemma}
\begin{proof}We claim that \[\div_{\TT\mathbb{S}^{2n-1}}\zeta=\div_{\R^{2n}}\zeta-\sum_{i=1}^{2n}
\left\langle\grad_{\R^{2n}}\zeta_i,
\frac{z}{\rho}\right\rangle\left\langle\frac{z}{\rho},
\ee_i\right\rangle,\]where $\zeta_i=\left\langle\zeta,
\ee_i\right\rangle$ for $i=1,..,2n.$ The last formula easily
follows by using definitions. Now the first term vanishes  because $\div_{\R^{2n}}z^\perp=0$. Using
\[\mathcal{J}_{\R^{2n}}\zeta=-\frac{1}{\rho}C\cc^{2n+1}-\frac{1}{\rho^4}\zeta\otimes z,\]yields
\[\sum_{i=1}^{2n}
\left\langle\grad_{\R^{2n}}\zeta_i,
\frac{z}{\rho}\right\rangle\left\langle\frac{z}{\rho},
\ee_i\right\rangle=\left\langle\mathcal{J}_{\R^{2n}}\zeta
\frac{z}{\rho}, \frac{z}{\rho}\right\rangle=0.\]
\end{proof}

\begin{lemma}\label{1lo}Under the previous assumptions, one has
\[\langle{\rm Hess}\cc\varphi\,  \nn,\nn\rangle=\rho^2\varphi_\rho''+(1-\rho^2)\varphi''_{\rho\zeta}+
\rho\sqrt{1-\rho^2}\varphi''_{\zeta \zeta}
+\frac{1-\rho^2}{\rho}\varphi'_\rho-\sqrt{1-\rho^2}\varphi'_\zeta.\]\end{lemma}

\begin{proof}By using Lemma \ref{tyie} we have
\begin{eqnarray}\label{1idrent}\langle{\rm Hess}\cc\varphi\,
\nn,\nn\rangle=\frac{\partial^2\varphi}{\partial\nn^2}+\left
\langle\frac{\qq\varpi}{\varpi},\qq
\varphi\right\rangle-\varpi\frac{\partial\varphi}{\partial\nn^\perp}.\end{eqnarray}We
compute \begin{eqnarray*}\frac{\partial^2\varphi}{\partial\nn^2}&=&
\frac{\partial}{\partial\nn}\left (\frac{\partial\varphi}{\partial
z} + \kappa\frac{\partial\varphi}{\partial
z^\perp}\right)\\&=&\frac{\partial}{\partial\nn}\left
(\rho\frac{\partial\varphi}{\partial \rho} +
 \sqrt{1-\rho^2}\frac{\partial\varphi}{\partial
\zeta} \right)\\&=&\frac{\partial}{\partial
z}\left(\rho\frac{\partial\varphi}{\partial \rho} +
 \sqrt{1-\rho^2}\frac{\partial\varphi}{\partial
\zeta}\right)+\kappa\frac{\partial}{\partial
z^\perp}\left(\rho\frac{\partial\varphi}{\partial \rho} +
 \sqrt{1-\rho^2}\frac{\partial\varphi}{\partial
\zeta}\right)\\&=&\rho\frac{\partial}{\partial
\rho}\left(\rho\frac{\partial\varphi}{\partial \rho} +
 \sqrt{1-\rho^2}\frac{\partial\varphi}{\partial
\zeta}\right)+ \sqrt{1-\rho^2}\frac{\partial}{\partial
\zeta}\left(\rho\frac{\partial\varphi}{\partial \rho} +
 \sqrt{1-\rho^2}\frac{\partial\varphi}{\partial
\zeta}\right).\end{eqnarray*} From this we easily get
that\begin{equation}\label{1pezzo}\frac{\partial^2\varphi}{\partial\nn^2}=\rho^2\varphi_{\rho
\rho}'' +
\rho\varphi'_\rho-\frac{\rho^2}{\sqrt{1-\rho^2}}\varphi'_\zeta +
2\rho\sqrt{1-\rho^2}\varphi''_{\zeta \rho} +
(1-\rho^2)\varphi''_{\zeta\zeta}.\end{equation}On the other hand,
using some computations  made in Section \ref{Sez3}, yields\[\frac{\qq\varpi}{\varpi}=-\frac{z\ss}{\rho^2(1-\rho^2)}.\]Moreover, one has
 \[z\ss=z-\langle z, \nn\rangle\nn=z-\rho^2\nn=z(1-\rho^2) -\rho\sqrt{1-\rho^2}z^\perp.\]
Therefore
\begin{eqnarray*}&&\left \langle\frac{\qq\varpi}{\varpi},\qq
\varphi\right\rangle-\varpi\frac{\partial\varphi}{\partial\nn^\perp}\\&=&-\left\langle\left(\frac{z(1-\rho^2)
-\rho\sqrt{1-\rho^2}z^\perp}{\rho^2(1-\rho^2)}+2\frac{\sqrt{1-\rho^2}}{\rho}\left(z^\perp-\kappa
z\right)\right),
\dg\varphi\right\rangle\\&=&-(2\rho^2-1)\left\langle \left(
\frac{z}{\rho^2}-\frac{z^\perp}{\rho\sqrt{1-\rho^2}}\right) ,
\dg\varphi\right\rangle\\&=&-(2\rho^2-1) \left(\frac{1}{\rho}
\frac{\partial\varphi}{\partial\rho}-\frac{1}{\sqrt{1-\rho^2}}\frac{\partial\varphi}{\partial\zeta}\right)
\\&=&\frac{1-2\rho^2}{\rho}\varphi'_\rho-\frac{1-2\rho^2}{\sqrt{1-\rho^2}}\varphi'_\zeta.\end{eqnarray*}
Using \eqref{1idrent}, \eqref{1pezzo} and the last computation
yields\begin{eqnarray*}&&\langle {\rm Hess}\cc
\varphi \,\nn,\nn\rangle\\&=&\rho^2\varphi_{\rho \rho}'' +
\rho\varphi'_\rho-\frac{\rho^2}{\sqrt{1-\rho^2}}\varphi'_\zeta + 2
\rho\sqrt{1-\rho^2}\varphi''_{\zeta \rho} +
(1-\rho^2)\varphi''_{\zeta \zeta}+
\frac{1-2\rho^2}{\rho}\varphi'_\rho-\frac{1-2\rho^2}{\sqrt{1-\rho^2}}\varphi'_\zeta
\\&=&\rho^2\varphi_\rho''+(1-\rho^2)\varphi''_{\zeta\zeta}+
2\rho\sqrt{1-\rho^2}\varphi''_{\zeta \rho}
+\frac{1-\rho^2}{\rho}\varphi'_\rho-\sqrt{1-\rho^2}\varphi'_\zeta \end{eqnarray*}which
achieves the proof.
\end{proof}

Remind that the Euclidean Laplace operator in
spherical coordinates is given by
$$\Delta_{\R^{2n}}\varphi
=\frac{1}{\rho^{2n-1}} \frac{\partial}{\partial
{\rho}}\left(\rho^{2n-1}\frac{\partial\varphi}{\partial
{\rho}}\right)
+\frac{1}{\rho^2}\Delta_{\mathbb{S}^{2n-1}}(\varphi|_{\mathbb{S}^{2n-1}(\rho)}),$$where $\Delta_{\mathbb{S}^{2n-1}}$ denotes the Laplacian on the Sphere ${\mathbb{S}^{2n-1}}$.
For the sake of simplicity, we will set
\[\Delta_{\mathbb{S}^{2n-1}}\varphi=\Delta_{\mathbb{S}^{2n-1}}(\varphi|_{\mathbb{S}^{2n-1}(\rho)}).\]

In particular, we have\[\Delta_{\R^{2n}}\varphi
=\varphi''_{\rho\rho}+\frac{2n-1}{\rho}\varphi'_\rho
+\frac{1}{\rho^2}\Delta_{\mathbb{S}^{2n-1}} \varphi.\]From
\eqref{1forma}, Lemma \ref{1lo} and the last formula we finally obtain\begin{eqnarray*}\lh \varphi&=& \Delta_{\R^{2n}}\varphi
-2\rho(n-\kappa^2)\varphi'_\rho-Q\rho\kappa\varphi'_\zeta-\left\langle{\rm
Hess}\cc\varphi \, \nn, \nn\right\rangle
\\&=&\varphi''_{\rho\rho}+\frac{2n-1}{\rho}\varphi'_\rho
+\frac{1}{\rho^2}\Delta_{\mathbb{S}^{2n-1}}
\varphi-2\rho(n-\kappa^2)\varphi'_\rho-Q\sqrt{1-\rho^2}\varphi'_\zeta\\&&-\left(\rho^2\varphi_\rho''+(1-\rho^2)\varphi''_{
\zeta\zeta}+ 2\rho\sqrt{1-\rho^2}\varphi''_{\zeta \rho}
+\frac{1-\rho^2}{\rho}\varphi'_\rho-\sqrt{1-\rho^2}\varphi'_\zeta\right)\\
&=& {(1-\rho^2)\varphi''_{\rho\rho}+\frac{2n
-(2n+1)\rho^2}{\rho}\varphi'_\rho}-
{2\rho\sqrt{1-\rho^2}\varphi''_{\zeta \rho}}\\&&+
{\frac{1}{\rho^2}\Delta_{\mathbb{S}^{2n-1}}\varphi
-(1-\rho^2)\varphi''_{\zeta\zeta}-(Q-1)
\sqrt{1-\rho^2}\varphi'_\zeta}.\end{eqnarray*}

We resume the previous computations in the next:
\begin{Prop}\label{ff} Let  $\varphi:\overline{B_1(0)}\setminus\{0\}\longrightarrow\R$ be a $\cont^2$-smooth function; then
the operator $\lh$ takes the following
form:\begin{eqnarray}\nonumber\lh\varphi&=&\underbrace{(1-\rho^2)\varphi''_{\rho\rho}+\frac{2n
-(2n+1)\rho^2}{\rho}\varphi'_\rho}_{\rm Radial\, Operator}-
\underbrace{2\rho\sqrt{1-\rho^2}\varphi''_{\zeta \rho}}_{\rm Mixed
\,Derivatives}\\\label{ffs}&&+\underbrace{\frac{1}{\rho^2}\Delta_{\mathbb{S}^{2n-1}}\varphi
-(1-\rho^2)\varphi''_{\zeta\zeta}-(Q-1)
\sqrt{1-\rho^2}\varphi'_\zeta}_{\rm
Angular\,Operator}.\end{eqnarray}Equivalently, one has
\begin{eqnarray}\nonumber\lh\varphi&=&\underbrace{(1-\rho^2)\left(\Delta_{\R^{2n}}\varphi-\varphi''_{\zeta\zeta}\right)
-2\rho\sqrt{1-\rho^2}\varphi''_{\zeta \rho}+
 \Delta_{\mathbb{S}^{2n-1}}\varphi}_{\rm 2nd\, order\, part}\\\label{ffso}&&
 \underbrace{-2\rho\varphi'_{\rho}-(Q-1)\sqrt{1-\rho^2}\varphi'_{\zeta}}_{\rm 1st\, order\, part}.\end{eqnarray}
\end{Prop}
\begin{proof}In order to prove formula \eqref{ffso}, it is enough using \eqref{ffs} together with the obvious identity
\[(1-\rho^2)\Delta_{\R^{2n}}\varphi=(1-\rho^2)\left(\varphi''_{\rho\rho}+\frac{2n-1}{\rho}\varphi'_\rho
+\frac{1}{\rho^2}\Delta_{\mathbb{S}^{2n-1}} \varphi\right).\]
\end{proof}
\begin{oss}[Case $\mathbb{H}^1$]If $n=1$,  fix polar coordinates on $\R^2$,
i.e. $(\rho, \vartheta)\in\R_+ \times[0, 2\pi]$ so that $z=(x, y)=(\rho\cos\vartheta, \rho\sin\vartheta)\in\R^2$. By applying Chain Rule, we get that
\[\frac{\partial\varphi}{\partial\zeta}=\frac{1}{\rho}\frac{\partial\varphi}{\partial\vartheta},\qquad
 \frac{\partial^2\varphi}{\partial\zeta^2}=\frac{1}{\rho^2}
\frac{\partial^2\varphi}{\partial\vartheta^2}\] for every smooth function
$\varphi:\R^2\longrightarrow\R$. We also stress that
\[\Delta_{\mathbb{S}^{1}}\varphi=\varphi''_{\vartheta\vartheta},\qquad  \frac{\partial^2\varphi}{\partial\zeta^2}=\frac{1}{\rho^2}
\frac{\partial^2\varphi}{\partial\vartheta^2}.\]By using
\eqref{ffso}, we compute
\begin{eqnarray}\nonumber \lh\varphi&=&(1-\rho^2)\Delta_{\R^{2}}\varphi
-2\rho\sqrt{1-\rho^2}\varphi''_{\zeta \rho}+
\rho^2\varphi''_{\zeta\zeta}-
 2\rho\varphi'_{\rho}-3\sqrt{1-\rho^2}\varphi'_{\zeta}\\&=&(1-\rho^2)\Delta_{\R^{2}}\varphi
-2\sqrt{1-\rho^2}\varphi''_{\vartheta \rho}+
\varphi''_{\vartheta\vartheta}-
 2\rho\varphi'_{\rho}-3\frac{\sqrt{1-\rho^2}}{\rho}\varphi'_{\vartheta}.\end{eqnarray}
Equivalently, by expressing the Laplacian $\Delta_{\R^{2}}$ in
polar coordinates, we see that
\begin{eqnarray}\label{aas} \lh\varphi&=&(1-\rho^2)\varphi''_{\rho\rho}- 2\sqrt{1-\rho^2}\varphi''_{\vartheta
\rho}+\varphi''_{\vartheta\vartheta} +\frac{2
-3\rho^2}{\rho}\varphi'_\rho-3
\frac{\sqrt{1-\rho^2}}{\rho}\varphi'_\vartheta.\end{eqnarray}
The operator $\lh$  turns out to be \textquotedblleft strictly non-elliptic\textquotedblright. In fact,  the determinant of the $(2\times 2)$-matrix
associated with the 2nd order coefficients  is zero at every point
of $\overline{B_1(0)}\setminus\{0\}$. In particular, by using
\eqref{1forma}, we easily see that the 2nd order part of
$\lh\varphi$ is given by\[\Delta_{\R^2}\varphi-\langle{\rm
Hess}_{\R^2} \varphi\, \nn,
\nn\rangle=(\nn^2)^2\frac{\partial^2\varphi}{\partial x^2} +
(\nn^1)^2\frac{\partial^2\varphi}{\partial
y^2}-2\nn^1\nn^2\frac{\partial^2\varphi}{\partial x \partial
y},\]where $\nn=(\nn^1, \nn^2)$. Since the matrix \[\left(%
\begin{array}{cc}
  (\nn^2)^2 & -\nn^1 \nn^2  \\
  -\nn^1 \nn^2  & (\nn^1)^2 \\
\end{array}%
\right),\]has eigenvalues $0$ and $1$ at every point of
$\overline{B_1(0)}\setminus\{0\}$, it follows that
$\lh$ is of \textquotedblleft parabolic-type\textquotedblright.
\end{oss}

Our final remarks  concern ``purely angular''
functions. So let
$\varphi:\overline{B_1(0)}\setminus\{0\}\longrightarrow\R$ be such that $\frac{\partial\varphi}{\partial\rho}=0.$
In case of purely angular functions,
Proposition \ref{ff} says that
\begin{eqnarray}\label{zaccazac}
\lh\varphi= {\frac{1}{\rho^2}\Delta_{\mathbb{S}^{2n-1}}\varphi
-(1-\rho^2)\varphi''_{\zeta\zeta}-(Q-1)
\sqrt{1-\rho^2}\varphi'_\zeta}.
\end{eqnarray}
Therefore, the eigenvalue equation becomes
\begin{eqnarray}\label{1spoint}{\frac{1}{\rho^2}\Delta_{\mathbb{S}^{2n-1}}\varphi
-(1-\rho^2)\varphi''_{\zeta\zeta}-(Q-1)
\sqrt{1-\rho^2}\varphi'_\zeta}=-\lambda\varphi.\end{eqnarray}
After multiplying both members of \eqref{1spoint} by $\rho^2$, we have\begin{eqnarray}\label{2spoint}
 \Delta_{\mathbb{S}^{2n-1}}\varphi-\rho^2(1-\rho^2)\varphi''_{\zeta\zeta}
-(Q-1)\rho^2\sqrt{1-\rho^2}\varphi_\zeta' =-\lambda\rho^2\varphi.
\end{eqnarray}

\begin{lemma}\label{mayerbeer}There exists no  non-trivial purely angular solution to Problem \ref{PEH}.\end{lemma}

\begin{proof}By contradiction. Changing the sign and differentiating  \eqref{2spoint} with respect to $\rho$, yields
\begin{eqnarray*}
 2\rho(1-2\rho^2)\varphi_\zeta''
+Q\frac{2-3\rho^2}{\kappa}\varphi_\zeta' =2\lambda\rho\varphi.
\end{eqnarray*}For a fixed $\rho\in]0, 1]$,  the last one is a 2nd order
O.D.E.  whose general solution must depend on $\rho$. This
contradicts the fact that $\varphi$ was assumed to be purely angular.\end{proof}

\begin{oss}Let $\varphi_i$ be the $i$-th radial eigenfunction of Problem
\ref{prad} associated with $\lambda_i$
 and let $\mu$ be an
angular function.  Note that
\[\lh(\varphi_i\mu)=\lh(\varphi_i)\mu+
\varphi_i\lh(\mu)+2\langle\qq\varphi_i,\qq\mu\rangle.\]Assume that $\varphi:=\varphi_i\,\mu$ is a solution of Problem
\ref{PEH} associated with  $\lambda$. Then, by applying
the previous formula we obtain\footnote{One uses also the
following:
\begin{eqnarray*}
2\langle\qq\varphi_i,\qq\mu\rangle=-2\frac{\partial\varphi_i}{\partial\nn}\frac{\partial\mu}{\partial\nn}=-2
\left\langle(\varphi_i)_{\rho}'\frac{z}{\rho},z+\kappa
z^\perp\right\rangle\langle\nabla_{\R^n}\mu,z+\kappa
z^\perp\rangle=-2\rho\sqrt{1-\rho^2}(\varphi_i)_{\rho}'{\mu}_\zeta'.
\end{eqnarray*}}
\begin{eqnarray*}\lambda_i\varphi_i\mu+\varphi_i\underbrace{\Big\{\frac{1}{\rho^2}
\Delta_{\mathbb{S}^{2n-1}}\mu-(1-\rho^2){\mu}_{\zeta\zeta}''
-(Q-1)\sqrt{1-\rho^2}{\mu}_\zeta'
\Big\}}_{=\lh\mu}-2\rho\sqrt{1-\rho^2}(\varphi_i)_{\rho}'{\mu}_\zeta'=-\lambda\varphi_i\mu.\end{eqnarray*}
So let $\rho\in]0,1[$ be such  that $\varphi_i(\rho)\neq 0.$ The previous equation gives
$$\lh\mu-2\frac{(\varphi_i)_{\rho}'\,\rho\,\sqrt{1-\rho^2}}{\varphi_i}{\mu}'_\zeta=-(\lambda-\lambda_i)\mu.$$
Multiplying by $-\mu$ and integrating the resulting equation along
$\Sph$ with respect to $\perh$, yields\begin{eqnarray*}\int_{\Sph}|\qq\mu|^2\,\perh
+\int_{\Sph}\left(\frac{(\varphi_i)_{\rho}'\,\rho\,\sqrt{1-\rho^2}}
{\varphi_i}\right)\frac{\partial(\mu^2)}{\partial\zeta}\,\perh=
(\lambda-\lambda_i)\int_{\Sph}\mu^2\,\perh.\end{eqnarray*}Furthermore, note that
$$\int_{\Sph}\left(\frac{(\varphi_i)_{\rho}'\,\rho\,\sqrt{1-\rho^2}}
{\varphi_i}\right)\frac{\partial(\mu^2)}{\partial\zeta}\,\perh=2\int_{0}^1\frac{(\varphi_i)_{\rho}'\sqrt{1-\rho^2}}
{\varphi_i}\frac{\rho^{2n+1}}{\sqrt{1-\rho^2}}\,d\rho\underbrace{\int_{\mathbb{S}^{2n-1}}
\frac{\partial(\mu^2)}{\partial\zeta}\,d\sigma_{\mathbb{S}^{2n-1}}}_{=0}=0,$$where
we have used the  Divergence Theorem for the Sphere
${\mathbb{S}^{2n-1}}$ together with Lemma \ref{rftr}. This implies that
$$\int_{\Sph}|\qq\mu|^2\,\perh=
(\lambda-\lambda_i)\int_{\Sph}\mu^2\,\perh,$$from which one obtains $\lambda\geq \lambda_i$. The last equation says
that $\mu$ must be an eigenfunction of $\lh$ on $\Sph$ with
eigenvalue $(\lambda-\lambda_i)$. By Lemma \ref{mayerbeer},  the only
possibility is $\lambda=\lambda_i$. It follows from Proposition
\ref{hopf} that $\mu$ must be constant.\end{oss}

The last remark shows how more general
solutions cannot be found in a too simple way. A possible strategy could be to use multiple Fourier series.\\

\begin{Prop}\label{teoauto}The spectrum of Problem \ref{prad} is contained in the spectrum of Problem \ref{PEH}.\end{Prop}

\begin{proof}Let $\varphi$ be a
smooth solution of Problem \ref{PEH} and fix spherical
coordinates $(\rho,\xi)\in[0, 1]\times\mathbb{S}^{2n-1}$. In
particular,  we are assuming that $\varphi$ solves the equation
$\lh\varphi=-\lambda\varphi$. Set
$\varphi_0(\rho):=\int_{\mathbb{S}^{2n-1}}\varphi(\rho,\xi)d\sigma_{\mathbb{S}^{2n-1}}(\xi).$
Roughly speaking, $\varphi_0$ is the spherical mean of $\varphi$ and hence a radial function.

At this point, the proof will immediately follow from the
next:\begin{claim}\label{yhyhyhyh} Let $\varphi_0\neq 0$. Then the function $\varphi_0$ is an
eigenfunction of Problem \ref{prad} associated with
$\lambda$.\end{claim}
\begin{proof}[Proof of Claim \ref{yhyhyhyh}]First, note that
$$\int_{\Sph}\varphi_0\,\perh=\int_{\Sph}\int_{\mathbb{S}^{2n-1}}\varphi\,\left(d\sigma_{\mathbb{S}^{2n-1}}\wedge\perh\right)=
\int_{\mathbb{S}^{2n-1}}\underbrace{\Big(\int_{\Sph}\varphi
\,\perh\Big)}_{=0}\,d\sigma_{\mathbb{S}^{2n-1}}=0.$$By making use
of \eqref{ffs} we have
\begin{eqnarray*}\lh\varphi&=& (1-\rho^2)\varphi''_{\rho\rho}+\frac{2n
-(2n+1)\rho^2}{\rho}\varphi'_\rho-
 2\rho\sqrt{1-\rho^2}\varphi''_{\zeta \rho} \\&&+\frac{1}{\rho^2}\Delta_{\mathbb{S}^{2n-1}}\varphi
-(1-\rho^2)\varphi''_{\zeta\zeta}-(Q-1)
\sqrt{1-\rho^2}\varphi'_\zeta.\end{eqnarray*}Integrating this
expression along $\mathbb{S}^{2n-1}$, using Lemma \ref{rftr} and
the Divergence Theorem for the Sphere ${\mathbb{S}^{2n-1}}$,
yields
\begin{equation}\label{hgty}\int_{\mathbb{S}^{2n-1}}\lh \varphi\,d\sigma_{\mathbb{S}^{2n-1}}=
\int_{\mathbb{S}^{2n-1}}\left((1-\rho^2)\varphi''_{\rho\rho}+\frac{2n
-(2n+1)\rho^2}{\rho}\varphi'_\rho\right)\,d\sigma_{\mathbb{S}^{2n-1}}.\end{equation}
Furthermore, one has
\begin{eqnarray*}\lh\varphi_0&=&(1-{\rho}^2)\frac{\partial^2\varphi_0}{\partial\rho^2}+\frac{2n
-(2n+1)\rho^2}{\rho}\frac{\partial\varphi_0}{\partial\rho}\\&=&\int_{\mathbb{S}^{2n-1}}\left((1-\rho^2)\varphi''_{\rho\rho}+\frac{2n
-(2n+1)\rho^2}{\rho}\varphi'_\rho\right)\,d\sigma_{\mathbb{S}^{2n-1}}
.\end{eqnarray*}So we finally get that
$$\lh\varphi_0=\int_{\mathbb{S}^{2n-1}}\lh \varphi
\,d\sigma_{\mathbb{S}^{2n-1}}=-\lambda\int_{\mathbb{S}^{2n-1}}\varphi\,d\sigma_{\mathbb{S}^{2n-1}}=-\lambda\varphi_0$$which
proves Claim \ref{yhyhyhyh}.\end{proof}
\end{proof}

\begin{oss}[Question]\label{zxc3} If $\mu$ denotes the 1st eigenvalue of Problem \ref{PEH}, then is it true that $\mu=Q-1$?
Roughly speaking, is the 1st eigenvalue of Problem \ref{PEH}  equal to the first eigenvalue of  Problem \ref{prad}?
\end{oss}

We now state a lemma
which is well-known in the classical setting; see \cite{FCS}.

\begin{lemma}\label{zxc}Let $S\subset\mathbb H^n$ be a hypersurface of class $\cont^2$ and let $\Om\subset S$ be any bounded open domain. If there exists a function $\psi>0$ on $\Om$ satisfying the equation $\lh \psi= q \psi$, then $$\int_{\Om}\left(|\qq \varphi|^2+ q \varphi^2\right)\,\perh\geq 0$$for all smooth function $\varphi$ compactly supported on $\Om$.
\end{lemma}

For a proof, see \cite{Montec}.

\begin{corollario}\label{zxc2}Let $\Om\Subset\Sph^+$ or $\Om\Subset\Sph^-$. Then, the following inequality holds
$$\int_{\Om}\left(|\qq \varphi|^2 - {(Q-1)}\,\varphi^2\right)\,\perh\geq 0$$for all smooth function $\varphi$ compactly supported on $\Om$.
\end{corollario}
\begin{proof} Setting $q=-(Q-1)$, the thesis follows by applying Lemma \ref{zxc} with the choice $\psi=\sqrt{1-\rho^2}$

\end{proof}

Another easy consequence of Lemma \ref{zxc} is contained in the next:
\begin{corollario}\label{zxc31}Let us set $\mathbb S^\ast:=\left\{p=\exp(z,t)\in\Sph: \rho=|z|\geq \sqrt{\frac{Q-1}{Q}}\right\}$. Then, for every
 $\Om\Subset\mathbb S^\ast$  the following inequality holds
$$\int_{\Om}\left(|\qq \varphi|^2 -  {2Q}\,\varphi^2\right)\,\perh\geq 0$$for all smooth function $\varphi$ compactly supported on $\Om$.
\end{corollario}
\begin{proof} Let $q=-{2Q}$. Note that the function $\psi:=\varphi_2(\rho)=Q\rho^2-(Q-1)$ is strictly positive on every open subset $\Om\Subset\mathbb S^\ast$. Since $\lh \psi= q \psi$, the thesis follows by applying Lemma \ref{zxc}.
\end{proof}

 In the inequalities of both Corollary \ref{zxc2} and Corollary \ref{zxc31}, the function $\varphi$ is not necessarily zero-mean. In other words, we do not require the validity of the condition $\int_\Om \varphi\,\perh=0$. Therefore, these inequalities are stronger than what one might expect: this seems to suggest a positive
answer to the question stated in Remark \ref{zxc3}. However, a
multiple Fourier series approach seems unavoidable in order to
give a rigorous proof.

\vspace{2cm}

{\footnotesize \noindent Francescopaolo Montefalcone:
\\Dipartimento di Matematica Pura e Applicata\\Universit\`a degli Studi di Padova,\\
  Via Trieste, 63, 35121 Padova (Italy)\,
 \\ {\it E-mail address}:  {\textsf montefal@math.unipd.it}}

\end{document}